\documentclass{amsproc}

\usepackage{amssymb}
\usepackage{graphicx}
\usepackage[cmtip,all]{xy}
\usepackage{amsmath}
\usepackage{amscd}
\usepackage{amsmath}

\usepackage{amsfonts,amssymb,amsopn}
\usepackage{amsthm}
\usepackage{mathrsfs}
\usepackage{color}
\usepackage{hyperref}
\usepackage{verbatim}

\newcommand{\scrD}{{\mathscr D}}
\newcommand{\cE}{{\mathcal E}}
\newcommand{\cF}{{\mathcal F}}

\newcommand{\cI}{{\mathcal I}}
\newcommand{\cK}{{\mathcal K}}
\newcommand{\cL}{{\mathcal L}}
\newcommand{\cO}{{\mathcal O}}
\newcommand{\cQ}{{\mathcal Q}}

\newcommand{\cT}{{\mathcal T}}
\newcommand{\cU}{{\mathcal U}}
\newcommand{\cV}{{\mathcal V}}
\newcommand{\cW}{{\mathcal W}}
\newcommand{\cX}{{\mathcal X}}
\newcommand{\cY}{{\mathcal Y}}
\newcommand{\cZ}{{\mathcal Z}}

\newcommand{\bA}{{\mathbb A}}
\newcommand{\C}{{\mathbb C}}

\newcommand{\PP}{{\mathbb P}}
\newcommand{\bQ}{{\mathbb Q}}
\newcommand{\bW}{{\mathbb W}}

\newcommand{\fm}{\mathfrak{m}}

\newcommand{\Sym}{\mathrm{Sym}}
\newcommand{\End}{\mathrm{End}\,}
\newcommand{\Ker}{\mathrm{Ker}}
\newcommand{\Aut}{\mathrm{Aut}\,}
\newcommand{\Image}{\mathrm{Im}}
\newcommand{\Iden}{\mathrm{Id}}
\newcommand{\Hom}{\mathrm{Hom}}
\newcommand{\rank}{\mathrm{rk}}
\newcommand{\Pic}{\mathrm{Pic}}

\newcommand{\isom}{\xrightarrow{\sim}}

\newcommand{\Kc}{{K_C}}
\newcommand{\Ox}{{\cO_X}}
\newcommand{\Oc}{{\cO_C}}

\newcommand{\Span}{\mathrm{Span}}
\newcommand{\defe}{\mathrm{def}}
\newcommand{\Spec}{\mathrm{Spec}\,}
\newcommand{\Quot}{\mathrm{Quot}}
\newcommand{\Hilb}{\mathrm{Hilb}}
\newcommand{\sm}{\mathrm{sm}}

\newcommand{\nd}{{\mathrm{nd}}}
\newcommand{\Supp}{\mathrm{Supp}}

\newcommand{\GL}{\mathrm{GL}}
\newcommand{\Gr}{\mathrm{Gr}}

\newcommand{\length}{\mathrm{length}}

\newcommand{\ope}{\cO_{\PP E}}
\newcommand{\opeo}{\ope (1)}
\newcommand{\Iz}{\cI_Z}
\newcommand{\Oz}{\cO_Z}
\newcommand{\Hef}{H^{e-f}_e}

\newcommand{\Hefl}{H^{e-f}_e ( \ell )}
\newcommand{\Heo}{H^{e-1}_e}
\newcommand{\Vef}{V^{e-f}_e}

\newcommand{\Vefl}{V^{e-f}_e ( \ell )}
\newcommand{\HeoV}{\Heo (V)}
\newcommand{\HerS}{\Hilb^e_\sm (S)}
\newcommand{\HefV}{\Hef (V)}
\newcommand{\Qef}{Q^{e-f}_e}
\newcommand{\QefV}{\Qef (V)}
\newcommand{\Qeo}{Q^{e-1}_e}
\newcommand{\QeoV}{\Qeo (V)}
\newcommand{\Ce}{C_e}

\newcommand{\hE}{{\widehat{E}}}
\newcommand{\hF}{{\widehat{F}}}
\newcommand{\hh}{{\widehat{h}}}

\newcommand{\htau}{{\widehat{\tau}}}
\newcommand{\tbW}{\widetilde{\bW}}

\newcommand{\tcX}{\widetilde{\cX}}
\newcommand{\qoe}{Q_{1,-e} (\Kc M^{-1} \otimes E)}
\newcommand{\inj}{\mathrm{inj}}
\newcommand{\Urde}{U_C ( r, d + e )}
\newcommand{\geM}{\gamma^e_M}
\newcommand{\hsigma}{\widehat{\sigma}}
\newcommand{\Me}{M^{[e]}}
\newcommand{\ut}{\underline{t}}
\newcommand{\uv}{\underline{v}}

\newtheorem{theorem}{Theorem}[section]
\newtheorem{lemma}[theorem]{Lemma}
\newtheorem{proposition}[theorem]{Proposition}
\newtheorem{corollary}[theorem]{Corollary}
\newtheorem{conjecture}[theorem]{Conjecture}

\theoremstyle{definition}
\newtheorem{definition}[theorem]{Definition}
\newtheorem{example}[theorem]{Example}

\theoremstyle{remark}
\newtheorem{remark}[theorem]{Remark}

\numberwithin{equation}{section}

\begin{document}

\title[Secant loci of scrolls]{Secant loci of scrolls over curves}


\author[G.\ H.\ Hitching]{George H.\ Hitching}
\address{Oslo Metropolitan University, Postboks 4, St. Olavs plass, 0130 Oslo, Norway}
\email{gehahi@oslomet.no}


\date{}

\begin{abstract}
Given a curve $C$ and a linear system $\ell$ on $C$, the secant locus $\Vef( \ell )$ parametrises effective divisors of degree $e$ which impose at most $e-f$ conditions on $\ell$. For $E \to C$ a vector bundle of rank $r$, we define determinantal subschemes $\Hefl \subseteq \Hilb^e ( \PP E )$ and $\QefV \subseteq \Quot^{0, e} ( E^* )$ which generalise $\Vef ( \ell )$, giving several examples. We describe the Zariski tangent spaces of $\QefV$, and give examples showing that smoothness of $\QefV$ is not necessarily controlled by injectivity of a Petri map. We generalise the Abel--Jacobi map and the notion of linear series to the context of Quot schemes.

We give some sufficient conditions for nonemptiness of generalised secant loci, and a criterion in the complete case when $f = 1$ in terms of the Segre invariant $s_1 (E)$. This leads to a geometric characterisation of semistability similar to that in \cite{Hit19}. Using these ideas, we also give a partial answer to a question of Lange on very ampleness of $\opeo$, and show that for any curve, $\QeoV$ is either empty or of the expected dimension for sufficiently general $E$ and $V$. When $\QeoV$ has and attains expected dimension zero, we use formulas of Oprea--Pandharipande and Stark to enumerate $\QeoV$.

We mention several possible avenues of further investigation.
\end{abstract}


\maketitle

\section{Introduction}

The purpose of this article is to investigate secant loci of scrolls over curves, and to begin a study of generalised secant loci on Quot schemes.

Let $C$ be a complex projective smooth curve of genus $g$. For $e \ge 1$, let $\Ce$ be the $e$th symmetric product of $C$. For $n \ge 1$, we consider the locus
\[ \Ce^f \ := \ \{ D \in \Ce : h^0 ( C, \Oc (D) ) \ge f + 1 \} \ = \ \{ D \in \Ce : h^0 ( C, \Kc (-D) ) \ge g - e + f \} . \]
This is a determinantal subvariety of $\Ce$, whose tangent spaces are defined by a Petri map. The construction and properties of $\Ce^f$ are explained in detail in \cite[Chapter IV]{ACGH}. The locus $\Ce^f$ is the preimage of the Brill--Noether locus
\[ W^f_e \ := \ \{ L \in \Pic^e (C) : h^0 ( C, L ) \ge f + 1 \} \]
by the Abel--Jacobi map $\Ce \to \Pic^e ( C )$.

An important characterisation of the divisors in $\Ce^f$ is given by the geometric Riemann--Roch theorem, which may be interpreted as saying that
\[ \Ce^f \ = \ \{ D \in \Ce : \dim \Span \left( \varphi_\Kc (D) \right) \le e - f - 1 \} . \]

More generally, let $L$ be a line bundle and $\ell = (L, V)$ a linear series of dimension $n$. For $D$ an effective divisor on $C$, we set
\[ \ell(-D) \ := \ \left( L, V \cap H^0 ( C, L(-D)) \right) . \]
For $0 \le f \le e$, we then have a \textsl{secant locus} associated to $\ell$, defined as
\begin{equation} \Vefl \ := \ \{ D \in \Ce : \dim \left( V \cap H^0 ( C, L ( -D )) \right) \ge n + 1 - e + f \} . \label{Classical} \end{equation}
Equivalently, writing $\varphi$ for the natural map $C \to \PP V^*$, we have
\[ \Vefl \ = \ \{ D \in \Ce : \dim \Span \left( \varphi (D) \right) \le e - f - 1 \} . \]
Note that $\Ce^f = \Vef \left( \Kc , H^0 ( C, \Kc ) \right)$.

The loci $\Ce^f$ are important invariants of the curve, as illustrated by the theorems of Martens and Mumford \cite[{\S} IV.5]{ACGH}, and their generalisations to $\Vefl$ proven in \cite{Baj15}. Various instances and properties of $\Ce^f$ have been investigated in \cite{Cop}, \cite{CM91}, \cite{GP82} and elsewhere. The tangent spaces of $\Ce^f$ and $\Vefl$ have been studied in \cite{CJ91}, \cite{CJ96}, \cite{Baj15} and \cite{Baj18}.

A major topic of interest is the enumerative geometry of $\Vefl$. This is studied for $\Ce^f$ in \cite[Chap.\ VII]{ACGH}, and that of $\Vefl$ is investigated in \cite{LeB06}, 
 \cite{Far08}, \cite{Far22}, \cite{Cot}, \cite{Ung21a}, \cite{Ung21b}, \cite{CHZ}, and elsewhere. In contrast to the Brill--Noether loci $W^f_e$, even when $\Vefl$ has nonnegative expected dimension it is not always defined by a map with suitable ampleness properties, so the question of nonemptiness may be nontrivial. In \cite{Ung19}, examples are given of empty secant loci with positive expected dimension, and in \cite{Ung21b} the question is addressed of when intersection-theoretic computations have enumerative meaning.

The purpose of the present work is to study two closely related generalisations of $\Vefl$, obtained by increasing $\dim (C)$ and $\rank (L)$ respectively. The first is very familiar. Let $S$ be a smooth projective variety of dimension $r \ge 1$, and $\ell = ( \cL, V )$ a linear series of dimension $n$ on $S$. Let $\HerS \subseteq \Hilb^e (S)$ be the component containing reduced subschemes. We define
\[ \Hefl \ := \ \{ Z \in \HerS : \dim \left( V \cap H^0 ( S, \cL \otimes \Iz ) \right) \ge n + 1 - e + f \} . \]
(The restriction to the smoothable component $\HerS$ is not necessary; see Remark \ref{restriction}.) Clearly, $\Hefl$ reduces to $\Vefl$ when $S$ is the curve $C$. In case $S$ is a surface, the enumerative geometry of $\Hefl$ has been studied via integrals of Segre classes in \cite{EGL}, \cite{Voi}, \cite{MOP} and elsewhere. Our focus is primarily on the case where $S$ is a projective bundle $\pi \colon \PP E \to C$ and $\cL = \opeo \otimes \pi^* M$, where techniques on vector bundles over curves can be brought to bear.

For the second generalisation of $\Vefl$, let $E \to C$ be a vector bundle of rank $r$ and degree $d$. The scheme $\Quot^{0, e} (E^*)$ parametrises torsion quotients of $E^*$ of length $e$; equivalently, subsheaves $F^* \subset E^*$ of rank $r$ and degree $-d - e$. Let $M$ be a line bundle, and suppose that $V \subseteq H^0 ( C, E^* \otimes M )$ is a subspace of dimension $n + 1$. We define
\[ \Qef ( E, M, V ) \ := \ \{ [ F^* \to E^* ] \in \Quot^{0, e} ( E^* ) : h^0 ( C, F^* \otimes M ) \ge n + 1 - e + f \} . \]
 If $E$ and $M$ are clear from the context, we write simply $\QefV$; and we abbreviate $\Qef \left( H^0 ( C, E^* \otimes M ) \right)$ to $\Qef$. Note that $\Quot^e ( L^{-1} ) = \Ce$ for all line bundles $L$. In this case, we have
\[ \Qef ( \Oc , L , V ) \ = \ \Qef ( L , \Oc , V ) \ = \ \Vef ( L^* , V ) . \]
Note in particular that
\begin{equation} \Ce^f \ = \ \Qef \left( \Kc , \Oc , H^0 ( C, \Kc ) \right) . \label{CefQK} \end{equation}

It is straightforward to see ({\S} \ref{defnHilb}) that $\Hefl$ and $\Qef (E, M, V)$, like $\Vefl$, are determinantal varieties with expected codimension $f ( n + 1 - e + f )$. Moreover, $\Qef ( H^0 ( C, E^* ) )$ is an analogue of a Brill--Noether locus for Quot schemes; see Remark \ref{BNQuot}. The loci $\QefV$ are better behaved than $\Hefl$ from many points of view. However, as the geometric intuition associated with $\Hefl$ is illuminating, we have as far as possible developed the two notions in parallel.\\
\par
Let us now give a summary of the article. After some preliminaries on scrolls and secant defect, we construct the loci $\Hefl$ and $\QefV$ and give some examples. In {\S} \ref{linkHQ}, using an approach from \cite{StaQS} and \cite{Hit20}, we describe a natural surjective rational map $\alpha \colon \Hef ( \opeo \otimes \pi^* M , V ) \dashrightarrow \Qef ( E, M, V )$ which is birational for $f = 0$. This is a useful tool in what follows. The indeterminacy locus of $\alpha$ is discussed briefly in {\S} \ref{pinondef}. We show in {\S} \ref{distinct} that $\QefV$ is not contained in $Q^{e-f-1}_e ( V )$ for $n + 1 - e + f \ge 0$, with a partial analogue for $\Hefl$. In {\S} \ref{GenProj} we discuss the effect of a general projection on secant loci both of type $\Hefl$ and $\QefV$.

Next, suppose that $\Quot^e ( E^* )$ contains a point $[ F^* \to E^* ]$ where $F$ is stable. Then there are rational classifying maps from $\HerS$ and $\Quot^e ( E^* )$ to the moduli space $U_C ( r, d+e)$ of stable bundles of rank $r$ and degree $d+e$. These generalise the Abel--Jacobi map $\Ce \to \Pic^e ( C )$, and map secant loci to Brill--Noether loci. In {\S} \ref{HefQefBN}, we study this situation and generalise the notion of linear series to the context of Quot schemes.

In {\S} \ref{nonempt}, we give some conditions for nonemptiness of the two types of generalised secant loci, extending various statements in the literature for $\Vefl$. In {\S} \ref{TangSp} the Zariski tangent spaces of $\Qef (E, M, V)$ are described. In the complete case, we show that smoothness is controlled by a Petri map; in general, however, the criterion for smoothness is not injectivity of the Petri map. Indeed, in {\S} \ref{examples} we show by example that when $\Qef$ admits a map to a Brill--Noether locus, smoothness of the secant locus is neither necessary nor sufficient for smoothness of the Brill--Noether locus in general.

For the remainder of the article, we focus on the case $f = 1$, generalising the study of inflectional loci of linearly normal scrolls in \cite{Hit19}; and using similar methods. In {\S} \ref{param}, we construct a parameter space for $\Qeo$ using Quot schemes of invertible subsheaves of $E$. This has several applications. In Proposition \ref{SegreInvChar}, we give a criterion for nonemptiness of $\left( \Heo \right)_\nd$ and $\Qeo$ in terms of the Segre invariant $s_1 (E)$. In Theorem \ref{LangeQ}, we offer a partial answer to a question of Lange \cite{Lan} on the relation between the Segre invariant and very ampleness of line bundles over ruled surfaces, which in fact is valid for projective bundles of any dimension over $C$. In {\S} \ref{SemistCohomSt}, we characterise semistable bundles in terms of secant loci. And in {\S} \ref{General}, we show that for a general bundle $E$ over any curve and a general $V \subseteq H^0 ( C, E^* \otimes M )$, the loci $\QefV$ are either empty or of the expected dimension. The approach for the latter is similar to \cite[{\S} 6]{Hit19}, and relies on Kleiman's theorem on transversality of translates. We conclude by showing how results from \cite{EGL}, \cite{OP} and \cite{StaCL} can be used to enumerate $\QeoV$ when it is has and attains expected dimension zero.

\subsection*{Topics of future inquiry}

Many questions which can be asked about Brill--Noether loci on $U_C ( r, d )$ or about moduli of coherent systems (see for example \cite{GT} and \cite{New}) can be formulated for $\Hefl$ and $\QefV$. For example, many questions of nonemptiness, components, dimension and smoothness of $\Hefl$ and $\QefV$ are open. The relation between smoothness of $\Qef$ and the Brill--Noether locus to which it maps, and the role of semistability in the whole, would be interesting to investigate further.

The enumerative geometry of $\QefV$ is of interest as a natural generalisation of that of $\Vefl$, and a natural application for the recent works \cite{OP} and \cite{StaCL}. It appears a priori more tractable than that of $\Hefl$: As a scroll contains many linear spaces, for $e \ge 3$ the loci $\Hef$ have components of excess dimension.

In another direction: Both $\Hefl$ and $\QefV$ appear naturally in the study of Mumford's notion of linear stability \cite{Mum} for scrolls in $\PP^n$. This has relevance for questions of moduli and, as explained in \cite{MS} and \cite{CT-L}, for Butler's well-known conjecture. This will be a topic of future investigation.

An object in some sense ``sitting between'' $\Hefl$ and the inflectional loci studied in \cite{Hit19} is the generalised Terracini locus defined in \cite{BC}. It would be interesting to see how the present techniques may be applicable to Terracini loci.

\subsection*{Acknowledgements} I thank the organisers of the VBAC 2022 meeting at the University of Warwick for financial support and for a very rewarding and enjoyable conference. I thank Ali Bajravani, Abel Castorena, Gavril Farkas and Peter Newstead for enjoyable and helpful communication. I thank Dragos Oprea for generous advice on enumeration questions. I am grateful to Samuel Stark for suggesting the strategy for Theorem \ref{counting} and for detailed explanation of the papers \cite{EGL}, \cite{OP} and \cite{StaCL}.

\section{Construction and first properties}

In this section we define generalised secant loci. After recalling or proving some background material on linear spans, secant defect and scrolls, we construct the loci $\Hefl$ and $\QefV$ and give some basic results. We give a number of examples.

\subsection{Linear spans and secant defect}

 Let $S$ be a variety and $\cL \to S$ a line bundle with nonempty linear system, and write $\varphi_\cL \colon S \dashrightarrow |\cL|^*$ for the standard map to the complete linear system of $\cL$. For any closed subscheme $Z \subseteq S$, let $r_Z \colon H^0 ( S, \cL ) \to H^0 ( Z, \cL|_Z )$ be the restriction map. Then it is easy to see that
\begin{multline} \Span \, \varphi_\cL ( Z) \ = \ \PP \Image \left( H^0 ( Z, \cL|_Z )^* \ \to \ H^0 ( S, \cL )^* \right) \ = \\
 \PP \Ker \left( H^0 ( S, \cL )^* \ \to \ H^0 ( S, \cI_Z \otimes \cL )^* \right). \label{SpanDescr} \end{multline}

Now let $V \subseteq H^0 ( S, \cL )$ be a nonzero subspace, and let $\psi = \colon S \dashrightarrow \PP V^*$ be the natural map. We have an exact diagram
\[ \xymatrix{ 0 \ar[r] & H^0 ( S, \cI_Z \otimes \cL ) \ar[r] & H^0 ( S, \cL ) \ar[r]^{r_Z} & H^0 ( Z , \cL|_Z ) \\
 0 \ar[r] & V \cap H^0 ( S, \cI_Z \otimes \cL ) \ar@{^{(}->}[u] \ar[r] & V \ar[r]^{r_Z} \ar@{^{(}->}[u] & r_Z ( V ) \ar@{^{(}->}[u] . } \]
Dualising and projectivising, we obtain a diagram of varieties
\begin{equation} \xymatrix{ \Span \, \varphi_\cL (Z) \ar@{^{(}->}[r] \ar@{-->}[d] & |\cL|^* \ar@{-->>}[r] \ar@{-->}[r] \ar@{-->}[d]^{p_V} & \PP H^0 ( S, \cI_Z \otimes \cL )^* \ar@{-->}[d] \\
 \Span \, \psi (Z) \ar@{^{(}->}[r] & \PP V^* \ar@{-->>}[r] \ar@{-->}[r] & \PP \left( V \cap H^0 ( S, \cI_Z \otimes \cL ) \right)^* . } \label{SpanDiag} \end{equation}
This description of $\Span \, \psi ( Z )$ will be useful in what follows.

We recall also the notion of secant defect of a zero-dimensional scheme.

\begin{definition} \label{DefectDefn} Let $Z \subset S$ be a subscheme of dimension zero. Let $\ell$ be a linear series on $S$ and $\psi \colon S \dashrightarrow \PP^n$ the associated map. The \textsl{defect} $\defe \, \psi ( Z )$ is the number
\[ \length \, Z - \dim \Span \, \psi (Z) - 1 . \]
(Here we take the empty set to have dimension $-1$.) We say that $Z$ is \textsl{$\ell$-nondefective} if $\defe \, \psi ( Z ) = 0$, and \textsl{$\ell$-defective} otherwise. \end{definition}

\subsection{Scrolls}

We briefly review scrolls over curves, which will be our primary objects of study. Let $C$ be a complex projective smooth curve of genus $g$. Let $E \to C$ be a vector bundle of rank $r$ and degree $d$, and $\pi \colon \PP E \to C$ the associated projective bundle. For $M \in \Pic (C)$, we write $\cL_M$ for the line bundle $\opeo \otimes \pi^* M$ over $\PP E$. By the projection formula and the definition of direct image, we have a natural identification $H^0 ( \PP E , \cL_M ) \isom H^0 ( C, E^* \otimes M )$. Fix a subspace $V \subseteq H^0 ( \PP E, \cL_M )$ of dimension $n + 1$, and consider the evaluation map
\begin{equation} \Oc \otimes V \ \to \ E^* \otimes M . \label{EvalMap} \end{equation}
Dualising, projectivising and projecting to $\PP V^*$, we obtain a map
\begin{equation} \psi \colon \PP E \ \dashrightarrow \ \PP V^* \times C \ \to \ \PP V^* . \label{DefnPsi} \end{equation}
This is naturally identified with the standard map $\PP E \dashrightarrow | \cL_M |^* \dashrightarrow \PP V^*$.

\subsection{Secant loci on Hilbert schemes} \label{defnHilb}

The first secant locus we will consider is a familiar and direct generalisation of the locus $\Ce^f \subset \Ce$. Let $S$ be a smooth projective variety of dimension $r \ge 1$. For $e \ge 1$, we write $\HerS$ for the component of $\Hilb^e (S)$ containing smoothable subschemes; that is, containing an open subset of $\Sym^e S$.

\begin{definition} Let $S$ be as above, and let $\ell = ( \cL , V )$ be a linear system of dimension $n$. Assume that $e \ge 1$ and $f \ge 0$ and $e \ge f$. We define
\[ \Hefl \ := \ \{ Z \in \HerS : \dim \left( \Span \, \psi_V ( Z ) \right) \le e - f - 1 \} \]
or, equivalently,
\begin{multline*} \Hefl \ = \ \{ Z \in \HerS : \defe \, \psi ( Z ) \ge f \} = \\
 \{ Z \in \HerS : \dim \left( V \cap H^0 ( S, \cL \otimes \cI_Z ) \right) \ge n + 1 - e + f \} . \end{multline*}
Depending on the context, we may write $\HefV$ instead. We write $\Hef$ for $\Hef \left( H^0 ( S, \cL) \right)$. \end{definition}

\begin{remark} \label{restriction} If $r \ge 3$ then $\Hilb^e (S)$ has multiple irreducible components. It is of course not necessary to restrict to the distinguished component $\HerS$; the study of secant loci on the other components of $\Hilb^e (S)$ is no less interesting. For the present, we restrict to $\HerS$ simply for convenience, as the other components are less well understood at the present time. \end{remark}

The locus $\Hefl$ can be constructed as a determinantal variety in a standard way. Let $p$ and $q$ be the projections of $\HerS \times S$ to the first and second factors respectively, and let $\cZ \subset \HerS \times S$ be the universal subscheme. Over $\HerS$ we have the diagram
\begin{equation} \xymatrix{ & \cO_{\HerS} \otimes V \ar[d] \ar[dr]^\varepsilon & \\
 p_* \left( \cI_\cZ \otimes q^* \cL \ar[r] \right) & p_* q^* \cL \ar[r] & p_* \left( \cO_\cZ \otimes q^* \cL \right) } . \label{HefDetStr} \end{equation}
As $\cZ$ is flat over $\Hilb^e ( S )$, the sheaf $p_* \left( \cO_\cZ \otimes q^* \cL \right)$ is locally free of rank $e$. Thus $\Hefl$ is the determinantal variety
\[ \{ Z \in \HerS : \rank ( \varepsilon|_Z ) \le e - f \} \]
defined by the $(e - f + 1)$st Fitting ideal of $\varepsilon$. The expected dimension of $\Hefl$ is therefore $re - f ( n + 1 - e + f )$.

\begin{example} \label{SumLineBundles} Let $L_1 , \ldots , L_r$ be line bundles on $C$. For $1 \le i \le r$, let $V_i \subseteq H^0 ( C, L_i )$ be a subspace of dimension $n_i + 1$. Set $E := \bigoplus_{i=1}^r L_i^{-1}$ and
\[ V \ := \ \bigoplus_{i=1}^r V_i \ \subseteq \ H^0 ( C, E^* ) \ = \ H^0 ( \PP E , \opeo ) , \]
and set $\ell = ( \opeo , V )$. Now for each $i$, we have a natural map
\[ \iota_i \colon C \ \to \ \PP V_i^* \ \hookrightarrow \ \PP \left( \bigoplus_{i=1}^r V_i^* \right) . \]
For each $i$, let $e_i, f_i$ be such that $V^{e_i - f_i}_{e_i} ( L_i , V_i )$ (defined in (\ref{Classical})) has a nonempty component of dimension $s_i$. Write $e := \sum e_i$ and $f := \sum f_i$. Then
\begin{equation} \left\{ \iota_1 ( D_1 ) \cup \cdots \cup \iota_r (D_r) : D_i \in V^{e_i - f_i}_{e_i} ( L_i, V_i ) \right\} \label{SumLBSecantLocus} \end{equation}
is a nonempty locus in $\Hefl$ of dimension $\sum s_i$. \end{example}

\noindent In {\S} \ref{pinondef}, Example \ref{SumLineBundles} will motivate a conjectural dimension bound for certain secant loci.

\begin{example} \label{Examplepidef} Using the fact that $\PP E$ is ruled, we can obtain loci in $\Hefl$ of excess dimension. Suppose $e \ge r + 1$. If $Z \in \HerS$ is reduced and has support in a fibre $\PP E|_x$ then $\Span \, \psi ( Z ) \subseteq \PP E|_x = \PP^{r-1}$. Hence any such $Z$ imposes at most $r$ conditions on sections of $\opeo$. Thus the locus
\[ \bigcup_{x \in C} \left\{ Z \text{ reduced of length } e : Z \subset \PP E|_x \right\} \]
belongs to $H^r_e (V)$ and has dimension
\[ \dim \left( \Sym^e \, \PP^{r-1} \right) + \dim C \ = \ e ( r - 1 ) + 1 \ = \ re - ( e - 1 ) . \]
This exceeds the expected dimension 
 $re - (e - r) ( n + 1 - r )$ of $H^r_e (V)$ if $2r \le n$. 
\end{example}

We will return in {\S} \ref{pinondef} to the study of subschemes with defect arising from points lining up inside a single fibre in this way.

\begin{remark} A linear series is base point free (resp., embedding) if and only if $H^0_1 ( \ell )$ (resp., $H^1_2 ( \ell )$) is empty. Most often we will consider embedded subvarieties of $\PP^n$, but many of our results apply to arbitrary projective models. \end{remark}

\subsection{Secant loci on Quot schemes}

Let us now give another generalisation of $\Vefl$. This time, instead of increasing $\dim (C)$, we increase $\rank (L)$.

\begin{definition} Let $E$ be a vector bundle of rank $r$ and degree $d$ over $C$, and let $\pi \colon \PP E \to C$ be the associated $\PP^{r-1}$-bundle. Let $V \subseteq H^0 ( C, E^* \otimes M )$ be a subspace of dimension $n + 1$. We define the generalised secant locus $\Qef ( E, M, V )$ as
\[ \left\{ [ F^* \to E^* ] \in \Quot^{0, e} ( E^* ) : \dim ( V \cap H^0 ( C, F^* \otimes M ) ) \ge n + 1 - e + f \right\} . \]
This can be constructed as a determinantal variety in a similar way to $\Hefl$. We may write simply $\Qef ( M, V )$ or $\QefV$ if no confusion should arise. Also, we often abbreviate $\Qef \left( H^0 ( C, E^* \otimes M ) \right)$ to $\Qef$. \end{definition}

\begin{remark} \label{BNQuot} 

As in \cite[{\S} 3.1]{HHN}, given a line bundle $M$ and any family of vector bundles $\cV \to B \times C$, one can define the \textsl{twisted Brill--Noether locus}
\[ B^k ( \cV , M ) \ = \ \{ b \in B : h^0 ( C, \cV_b \otimes M ) \ge k \} \ \subseteq \ B . \]
Let $\cF^*$ be the universal subsheaf over $\Quot^{0, e} (E^*) \times C$. Then we have simply
\[ \Qef \left( E, M, H^0 ( C, E^* \otimes M ) \right) \ = \ \Qef \ = \ B^{h^0 ( C, E^* \otimes M )-e+f} ( \cF^* , M ) . \]
Other relations between secant loci and Brill--Noether loci will be studied in {\S} \ref{HefQefBN}.
\end{remark}

\subsection{The link between \texorpdfstring{$\HefV$}{H e-f e} and \texorpdfstring{$\QefV$}{Q e-f e}} \label{linkHQ}

Following \cite{Hit20} and \cite{StaQS}, for $S = \PP E$ we now describe the relation between the two generalised secant loci above. This will generalise the fact that if $L$ is a line bundle, then
\[ D \mapsto \left[ L^{-1} (-D) \to L^{-1} \right] \text{ defines an isomorphism } \Ce \isom \Quot^{0, e} ( L^{-1} ) . \]
Let $Z \subset \PP E$ be a subscheme of length $e$. Taking direct images of $0 \to \Iz (1) \to \opeo \to \Oz (1) \to 0$ on $C$, we have a sequence
\[ 0 \ \to \ \pi_* \Iz (1) \ \to \ E^* \ \to \ \pi_* \Oz (1) \ \to \ \cdots \label{DIseq} \]
As $\pi_* \Oz (1)$ is a torsion sheaf on $C$, we see that $\pi_* \Iz (1)$
is a full rank subsheaf of $E^*$. Following \cite{Hit20}, we denote this by $E_Z^*$. Dualising, for general $Z \in \HerS$, we have $\deg (E_Z) = \deg (E) + e$; for example if $\pi(Z)$ consists of $e$ distinct points.

\begin{proposition} \label{alphaProperties} Let $E \to C$ be a vector bundle.
\begin{enumerate}
\item[(a)] The association $Z \mapsto \left[ E_Z^* \to E^* \right]$ defines a rational map $\alpha \colon \Hilb^e (S) \dashrightarrow \Quot^{0, e} (E^*)$.
\item[(b)] The restriction of $\alpha$ to the component $\HerS$ is surjective. 
 Moreover, $\alpha$ is bijective on the open set
\[ \{ Z \in \HerS : \Supp \, \pi (D) \hbox{ consists of $e$ distinct points} \} . \]
In particular, $\alpha|_{\HerS}$ is a birational equivalence.
\end{enumerate}
\end{proposition}

\begin{proof} (a) The existence of a suitable family of length $e$ quotients $E^* \to \pi_* \Oz ( 1 )$ parametrised by $\Hilb^e ( \PP E )$ is proven\footnote{It was stated, unfortunately without an adequate justification, in \cite[Theorem 2.6 (a)]{Hit20}.} in \cite[Theorem 3 (i)]{StaQS} for $e = 2$ when $\PP E$ is a projective bundle over a surface, and the same argument works for $e \ge 1$ when the base is a curve. In this case, over the open subset
\[ \{ Z \in \Hilb^e (S) : E^* \to \pi_* \Oz (1) \hbox{ is surjective} \} \]
we obtain a family of elementary transformations $E_Z^* \to E^*$ of degree $\deg (E^*) - e$. This locus is nonempty, as it contains all $Z$ projecting to $e$ distinct points of $C$, so we obtain the desired rational map $\alpha$.

(b) This follows from the proof of \cite[Theorem 2.6]{Hit20}. \end{proof}

\begin{remark} If $e \ge 2$ and $r \ge 2$, then $\alpha$ has fibres of positive dimension. For example, suppose $e = 2$ and let $\mu_1$ and $\mu_2$ be points of a fibre $\PP E|_x$ spanning a $\PP^1$. Then for any $\mu_1', \mu_2'$ spanning the same $\PP^1$, we have $E_{\{ \mu_1' , \mu_2' \}} = E_{\{ \mu_1 , \mu_2 \}}$. \end{remark}

\subsubsection{Relatively defective subschemes} \label{pinondef}

In contrast to the situation studied in \cite{StaQS}, the map $\alpha$ is not a morphism when $e \ge 3$. The indeterminacy locus is the set of $Z \in \Hilb^e (S)$ such that $E^* \to \pi_* ( \Oz \otimes \opeo )$ is not surjective; equivalently, $\deg ( E_Z^* ) > \deg E^* - e$. For example, if $Z$ consists of three collinear points in a fibre of $\PP E$, then $\deg (E_Z^*) = \deg E^* - 2$. This motivates the following, which is \cite[Definition 3.4]{Hit20}:

\begin{definition} \label{DefnPiDef} Let $Z \subset \PP E$ be a subscheme of dimension zero. Then $Z$ is said to be \textsl{$\pi$-nondefective} if $\deg ( E_Z^* ) = \deg ( E^* ) - \length (Z)$, and \textsl{$\pi$-nondefective} otherwise. \end{definition}

Now we can formulate an important statement and a conjecture. We denote by $\HerS_\nd$ the locus of $\pi$-nondefective subschemes in $\HerS$, and we write
\[ \HefV_\nd \ := \ \HefV \cap {\HerS}_\nd. \]

\begin{proposition} Let $V$ be an $n+1$-dimensional subspace of
\[ H^0 ( C, E^* \otimes M ) \ \cong \ H^0 ( S, \opeo \otimes \pi^* M ) . \]
Suppose $n + 1 \ge e - f$. Let $Z \in \HerS$ be $\pi$-nondefective. Then $Z \in \HefV$ if and only if $[E_Z^* \to E^*]$ belongs to $\QefV$. \label{HQequiv} \end{proposition}

\begin{proof} By construction of $E_Z$, we have
\begin{equation} H^0 ( C, E_Z^* \otimes M ) \cap V \ = \ H^0 ( S , \cI_Z (1) \otimes \pi^* M ) \cap V . \label{EZIZ} \end{equation}
The statement follows. \end{proof}

In view of the proposition, we will switch freely between $\HefV_\nd$ and $\QefV$, depending on what is more convenient or illuminating. We refer to \cite[{\S} 3.1]{Hit20} for discussion and some equivalent definitions of $\pi$-nondefectivity. Example \ref{Examplepidef} illustrates that $\pi$-defective subschemes can give rise to artificially large components of $\Hefl$, which do not reflect the behaviour of $\QefV$. For the most part, we will exclude $\pi$-defective subschemes from our consideration.

Next, we turn briefly to the question of dimension bounds for $\HefV$ and $\QefV$. Returning to Example \ref{SumLineBundles}: If for each $i$ the pair $(C, L_i)$ is one of the four special types described in \cite[Theorem 4.6]{Baj15} and $e_i \le h^0 ( C, L_i ) - 2$, then we may take $s_i$ to be the maximal value $e_i - f_i - 1$ proven in loc.\ cit., and then the locus (\ref{SumLBSecantLocus}) above has dimension $e - f - r$. Based on this information, we make the following naive conjecture.

\begin{conjecture} Suppose $C$ is nonhyperelliptic of genus $g \ge 9$. Let $E \to C$ be a bundle of rank $r$ and slope at most $g-2$ such that $\opeo$ is very ample. For $e \le h^0 ( C, E ) - 2r$, we have
\[ \dim \Qef \ \le \ \dim \left( \Hef \right)_\nd \ \le \ e - f - r . \]
\end{conjecture}

We also remark the following fact, which generalises \cite[Proposition 4.4 and Remark 4.5 (a)]{Hit20} in the complete case.

\begin{lemma} 
 Let $[ F^* \to E^* ]$ be a point of $\Quot^{0, e} (E^*)$. Let $Z \subseteq \PP E$ be any subscheme such that $F = E_Z$. Then
\[ \Span \, \psi ( Z ) \ = \ \PP \Ker \left( V^* \ \to \ ( V \cap H^0 ( C, E_Z^* \otimes M )^* \right) \ = \ \PP \left( V \cap H^0 ( C, E_Z^* \otimes M ) \right)^\perp . \]
In particular, $\Span \, \psi (Z)$ depends only on $F^*$, and in fact does not require that $Z$ be $\pi$-nondefective or even zero-dimensional
\end{lemma}

\begin{proof} The statement is a consequence of (\ref{SpanDiag}) and (\ref{EZIZ}). \end{proof}

\subsection{Distinctness of the strata} \label{distinct}

For fixed $e$, the determinantal varieties $\Hefl$ and $\QefV$ define stratifications of $\HerS$ and $\Quot^{0, e} (E^*)$ respectively as $f$ varies. Generalising \cite[Lemma 4.1]{Baj15}, \cite[Lemma 2.2]{AS} and \cite[Lemma IV.1.7]{ACGH}, we will now show that $\QefV$ is strictly bigger than $Q^{e-f-1}_e (V)$ whenever $n + 1 - e + f \ge 0$, and similarly for $\Hefl$ under a certain condition. The result for higher dimensional varieties rests on the following statement for curves, which is a straightforward generalisation of \cite[Lemma 2.2]{AS}; as the proof is essentially identical, we omit it.

\begin{lemma} \label{firstpropscurve} Let $T$ be a smooth projective curve and $\ell = ( \cL , V )$ a very ample linear system on $T$ of dimension $n + 1$. Let $e$ and $f$ be integers satisfying $0 \le f \le e$ and $n + 1 - e + f \ge 0$. Then no component of $\Vefl$ is fully contained in $V^{e-f-1}_e (\ell)$.\end{lemma}

\begin{proposition} \label{firstprops} Let $S$ be a smooth projective variety, and $\ell = ( \cL , V )$ a very ample linear system of dimension $n + 1$. Let $e$ and $f$ be integers with $e \ge 1$ and $f \ge 0$ and $e \ge f$, and satisfying $n + 1 - e + f \ge 0$.

\begin{enumerate}
\item[(a)] Let $\cX$ be a component of $\Hefl$. Suppose $Z_0 \in \cX$ is such that the subscheme $Z_0 \subset S$ is curvilinear. Then $\cX$ is not fully contained in $H^{e-f-1}_e ( \ell )$.
\item[(b)] Let $E \to C$ be a vector bundle, and $\cL_M = \opeo \otimes \pi^* M$. Let $V \subseteq H^0 ( C, E^* \otimes M )$ be a subspace of dimension $n + 1$. Then no component of $\QefV$ is fully contained in $Q^{e-f-1}_e (V)$.
\end{enumerate}
\end{proposition}

\begin{proof} (a) As curvilinearity is an open condition, clearly we may assume that $Z_0$ belongs to no other component of $\Hefl$ than $\cX$. 
 By curvilinearity, $T_\nu Z_0$ has dimension one for all $\nu \in \Supp (Z_0)$. Hence, by considering the closure of a suitable affine curve contained in some affine open subset of $S$, we may assume that $Z_0$ belongs to $\Hilb^e (T) \cong T_e$ for a curve $T \subseteq S$ whose image in $\PP V^*$ is nondegenerate. Then we obtain an identification of $V$ with a subspace of $H^0 \left( T, \cL|_T \right)$. Setting $\ell_T := \left( \cL|_T, V \right)$, we can consider the usual secant locus $V^{e-f}_e ( \ell_T )$ as defined in (\ref{Classical}). For all $f \ge 0$, there are natural inclusions $V^{e-f}_e ( \ell_T ) \hookrightarrow \Hefl$. As $Z_0$ belongs only to the component $\cX$ of $\Hefl$, if $\cV$ is any component of $V^{e-f}_e ( \ell_T )$ containing $D_0$ then we have $\cV \subseteq \cX$.

Now suppose that $Z_0 \in H^{e-f-1}_e ( \ell )$. Then $Z_0 \in V^{e-f-1}_e ( \ell_T ) \cap \cV$ also. Since $e \le n + 1 + f$, by Lemma \ref{firstpropscurve} we may deform $Z_0$ inside $\cV$ to a point $Z_1 \in \cV \setminus V^{e-f-1}_e ( \ell )$. Then $Z_1$ has defect precisely $f$ in $\PP V^*$. But as $\cV \subseteq \cX$, we see that $\cX$ contains the point $Z_1 \in \Hefl \setminus H^{e-f-1}_e ( \ell )$. This proves (a).

(b) Let $\cX$ be a component of $\QefV$. By Lemma \ref{alphaProperties} (a) and Proposition \ref{HQequiv}, we see that $\alpha^{-1} ( \cX )$ is a union $\bigcup_i \tcX_i$ of nonempty open subsets of components of $\HefV_\nd$. By the proof of \cite[Theorem 2.6 (a)]{Hit20}, at least one $\tcX_{i_0}$ contains a point $Z_0$ corresponding to a curvilinear subscheme. By part (a), therefore, $\tcX_{i_0}$ is not contained in $H^{e-f-1}_e (V)$. Again by Lemma \ref{HQequiv}, the image of any point of $\tcX_{i_0} \setminus H^{e-f-1}_e (V)$ is a point of $\cX \setminus Q^{e-f-1}_e (V)$, and we obtain (b). \end{proof}

\subsection{Secant loci under general projections} \label{GenProj}

We will now show that the secant loci $\HefV$ and $\QefV$ behave predictably under projection from a general centre. The main ingredient is the theorem of Kleiman on transversality of a general translate, which will be used several times in what follows:

\begin{theorem}[{\cite[Theorem 2 (i)]{Kle}}] Let $G$ be a connected algebraic group and $Z$ an irreducible variety admitting a transitive action of $G$. Suppose that $a \colon X \to Z$ and $b \colon Y \to Z$ are maps of irreducible varieties. For any $s \in G$, we write $\gamma \cdot Y$ for $Y$ considered as a $Z$-variety via the map $y \mapsto \gamma \cdot b ( y )$. Then there is a dense open subset $U \subseteq G$ such that $( \gamma \cdot Y ) \times_Z X$ is either empty or equidimensional of dimension $\dim X + \dim Y - \dim Z$ for all $\gamma \in U$. \label{KleimanTheorem} \end{theorem}

Our results will follow from this theorem and linear algebra. To avoid repetition, \textbf{during this subsection, we fix integers $n$, $m$, $e$, and $f$ such that}
\[ 0 \ \le \ f \ \le \ e \ \le \ m + 1 \ \le \ n + 1 . \]

Let $H$ be a variety and $\beta \colon \cE \to \cO_H \otimes V^*$ a map of vector bundles over $H$, where $\rank \, E = e$ and $V^* = \C^{n+1}$. We consider the degeneracy loci
\[ D_{e-f} ( \beta ) \ = \ \{ h \in H : \rank \left( \beta|_h \right) \le e - f \} \quad \hbox{and} \quad D_{e-f}^\circ ( \beta ) \ = \ \{ h \in H : \rank \left( \beta|_h \right) = e - f \} . \]
For each $f'$ with $0 \le f' \le f$, we write
\[ d_1 ( f' ) \ := \ \dim B - f' ( n + 1 - e + f' ) \quad \hbox{and} \quad d_2 ( f' ) \ := \ \dim B - f' ( m + 1 - e + f ) . \]
Notice that $d_1 (f')$ is the expected dimension of $D_{e-f'} ( \beta )$ and $D_{e-f'}^\circ ( \beta )$. Lastly, if $W \subseteq V$ is a subspace, abusing notation, we write $p_W$ for the projections $V^* \to W^*$ and $\cO_H \otimes V^* \to \cO_H \otimes W^*$, and $W^\perp := \Ker ( p_W )$.

\begin{proposition} \label{ProjOfDetVar} Let $H$ and $\beta \colon \cE \to \cO_H \otimes V^*$ be as above. Assume that for $0 \le f' \le f$, we have
\begin{equation} \dim D_{e-f'} ( \beta ) \ \le \ \max\{ d_1 ( f' ) , d_2 ( f' ) \} . \label{technical} \end{equation}
Then for a general subspace $W \subseteq V$ of dimension $m + 1$, the degeneracy locus $D_{e-f} ( p_W \circ \beta )$ is empty or of the expected dimension $\dim H - f ( m + 1 - e + f )$ when this is nonnegative. \end{proposition}

\begin{proof} Suppose $h \in D_{e-f'}^\circ ( \beta )$ for some $f' \in \{ 0 , \ldots , f \}$. Then $\rank ( p_W \circ \beta|_h ) \le e - f$ if and only if
\begin{equation} \dim \left( W^\perp \cap \Image \, \beta|_h \right) \ \ge \ f - f' . \label{SpanInters} \end{equation}
Firstly, suppose $f - f' > n - m = \dim W^\perp$. Then condition (\ref{SpanInters}) is not satisfied for any $h \in D_{e - f'}^\circ ( \beta )$, so $D_{e-f} (p_W \circ \beta|_h) \cap D_{e - f'}^\circ ( \beta )$ is empty for all $W$ of dimension $m + 1$.

 On the other hand, suppose that $f - f' \le n - m$. (One can show that this implies that $d_2 ( f' ) \ge d_1 ( f' )$.)
 Let $a \colon D_{e-f'}^\circ ( \beta ) \to \Gr ( e - f' , V^* )$ be the map $h \mapsto \Image \, \beta|_h$. Consider the Schubert cycle
\[ \Sigma^{f-f'}_{W^\perp} \ := \ \{ \Lambda \in \Gr ( e-f' , V^* ) : \dim ( \Lambda \cap W^\perp ) \ge f - f' \} . \]
In the notation of \cite[p.\ 194--196]{GH}, we have $\Sigma^{f-f'}_{W^\perp} = \overline{W_{a_1 , \ldots , a_{f-f'}}}$ where $a_1 = \cdots = a_{f-f'} = m + 1 - e + f$. By loc.\ cit.,
\begin{equation} \dim \Sigma^{f-f'}_{W^\perp} \ = \ \dim \Gr ( e - f' , V^* ) - (f-f') ( m + 1 - e + f ) . \label{dimSigma} \end{equation}
 For each component $Y$ of $D_{e-f'}^\circ ( \beta )$, we consider the fibre product diagram
\[ \xymatrix{ Y_W \ar[r] \ar[d] & Y \ar[d]^{a|_Y} \\
 \Sigma^{f-f'}_{W^\perp} \ar@{^{(}->}[r] & \Gr ( e - f' , V^* ) . } \]
If $h \in Y$, then by (\ref{SpanInters}) we see that $h \in D_{e-f} ( p_W \circ \beta )$ if and only if $h \in Y_W$.

Now $\GL ( V^* )$ acts transitively on $\Gr ( \ell , V^* )$ for all $\ell$. For $\gamma \in \GL ( V^* )$, unwinding definitions, we see that $\gamma \cdot \Sigma^{f-f'}_{W^\perp} \ = \ \Sigma^{f-f'}_{\gamma \cdot W^\perp}$. Hence by Theorem \ref{KleimanTheorem} and by (\ref{dimSigma}), for a general choice of $W$ the locus $Y_W$ is empty or equidimensional of dimension
\[ \dim Y - (f - f') ( m + 1 - e + f ) . \]
As by hypothesis (\ref{technical}) we have $\dim Y \le d_2 ( f' )$, this is bounded above by
 the expected dimension of $D_{e-f} ( p_W \circ \beta )$. (In fact, as the latter is determinantal, $Y_W$ is nonempty only if $\dim Y = d_2 ( f' )$.) \end{proof}

This has the following applications to the study of secant loci. The following generalises \cite[Theorem A.1]{Hit19}. It may be viewed as an elementary contribution to the line of inquiry in for example \cite{GP13} and \cite{Ran}. For a projective model $\psi \colon S \dashrightarrow \PP V^*$ and $W \subseteq V$, let $\psi_W$ be the composition $S \dashrightarrow \PP V^* \dashrightarrow \PP W^*$.

\begin{corollary} Let $S$ be a projective variety, and $\psi \colon S \dashrightarrow \PP^n = \PP V^*$ a projective model. Let $\cZ \subset H \times S$ be a family of length $e$ subschemes of $S$ parametrised by a variety $H$. For each subspace $W \subseteq V$ and for each $f'$, we define
\[ H^{e-f'}_e ( H , W ) \ := \ \{ h \in H : \defe \, \psi_W ( \cZ_h ) \ge f' \} . \]
Assume that for $0 \le f' \le f$ we have $\dim H^{e-f'}_e ( H, V ) \le \max \{ d_1 ( f' ) , d_2 ( f' ) \}$. Then for a general subspace $W \subseteq V$ of dimension $m + 1$, the secant locus $\Hef ( H , W )$ is empty or of the expected dimension $\dim H - f ( m + 1 - e + f )$. \end{corollary}

\begin{proof} For any $W \subseteq V$, the locus $H^{e-f'}_e (H, W)$ can be constructed as a determinantal variety exactly as in {\S} \ref{defnHilb}. Let $\cL$ be the line bundle inducing the map to $\PP^n$. Let $p$ and $q$ be the projections of $H \times S$ to the first and second factors respectively. We have a diagram of sheaves over $H \times S$ as follows:
\[ \xymatrix{ & \cO_H \otimes W \ar[d]_{^tp_W} \ar[ddr] & \\
 & \cO_H \otimes V \ar[d] \ar[dr]_\varepsilon & \\
 p_* \left( \cI_\cZ \otimes q^* \cL \right) \ar[r] & p_* q^* \cL \ar[r] & p_* \left( \cO_\cZ \otimes q^* \cL \right) } . \]
As $\cZ$ is flat over $H$, the sheaf $p_* \left( \cO_\cZ \otimes q^* \cL \right)$ is locally free of rank $e$. Now for $h \in H$, by (\ref{SpanDescr}) we have $\Span \, \psi_W ( \cZ_h ) = \PP \Image \left( p_W \circ {^t\varepsilon|_h} \right)$. Thus $H^{e-f'}_e (H, W)$ is the determinantal variety
\[ \{ Z \in H : \rank ( p_W \circ {^t\varepsilon} )|_Z \le e - f' \} \ = \ D_{e-f'} ( p_W \circ {^t\varepsilon} ) . \]

Now by hypothesis, condition (\ref{technical}) is satisfied by $D_{e-f'} ( {^t\varepsilon} )$ for $0 \le f' \le f$. Therefore, $\Hef ( H, W ) = D_{e-f} ( p_W \circ {^t\varepsilon} )$ is empty or of the expected dimension by Proposition \ref{ProjOfDetVar} when this is nonnegative, and empty otherwise. \end{proof}

We will now give a corresponding statement for secant loci of type $\QefV$. This will be used in {\S} \ref{General}.

\begin{corollary} \label{QefGenProj} Let $C$ be a curve and $E \to C$ a vector bundle. Let $V \subseteq H^0 ( C, E^* )$ be a subspace of dimension $n+1$. Let $\cF^\vee \to p_C^* E^*$ be a family of elementary transformations of degree $\deg E^* - e$ parametrised by a variety $Q$. For each $f'$ and for each subspace $W \subseteq V$ of dimension $m + 1$, set
\[ Q^{e-f'}_e ( Q , W ) \ := \ \{ q \in Q : h^0 ( C, \cF^*_q ) \cap W \ge m + 1 - e + f' \} . \]
Assume that for $0 \le f' \le f$ we have $\dim Q^{e-f'}_e ( Q , V ) \le \max \{ d_1 ( f' ) , d_2 ( f' ) \}$. Then for a general subspace $W \subseteq V$ of dimension $m + 1$, the locus $\Qef ( Q , W )$ is empty or of the expected dimension $\dim Q - f ( m + 1 - e + f )$ when this is nonnegative, and empty otherwise. \end{corollary}

\begin{proof} This is similar to the proof of the previous corollary. \end{proof}

This completes our overview of basic properties of scrolls and their secant loci. In the sections that follow, we will apply these to various questions on $\Hefl$ and $\QefV$.

\section{Secant loci and Brill--Noether loci} \label{HefQefBN}

We now further discuss the connection between secant loci and Brill--Noether loci, and generalise the notions of Abel--Jacobi map and linear series to the context of Quot schemes. This expands \cite[Remark 4.7]{Hit20}.

It is well known that for $e \ge f \ge 0$ there is a diagram
\begin{equation} \xymatrix{ \Ce \ar[r]^-\sim & \Quot^{0, e} ( \Oc ) \ar[r] & \Pic^e ( C ) \\
 \Ce^f \ar[r]^-\sim \ar@{^{(}->}[u] & \Qef \left( \Oc , \Kc , H^0 ( C, \Kc ) \right) \ar[r] \ar@{^{(}->}[u] & B_{1, e}^{1 + f} \ar@{^{(}->}[u] } \label{LineBundleHefQef} \end{equation}
The rows are given by $D \mapsto [ \Oc (-D) \to \Oc ] \mapsto \Oc ( D )$. This composition is an Abel--Jacobi map, with fibre over $L$ given by the linear series $|L|$.

Now let $\Urde$ be the moduli space of stable bundles of rank $r$ and degree $d+e$ over $C$. For $k \ge 0$, we consider the higher rank Brill--Noether locus
\[ B_{r, d+e}^k \ = \ \{ F \in \Urde : h^0 ( F ) \ge k \} . \]
Let $E$ be a vector bundle of rank $r$ and degree $d$, not necessarily semistable. Suppose that $0 \le f \le e$ and $\Qef ( E, \Kc , H^0 ( E^* \otimes \Kc ) )$ contains a point $[ F^* \to E^*]$ such that $F$ is a stable bundle. Now by Serre duality, there is an exact sequence of vector spaces
\begin{multline*} 0 \ \to \ H^0 ( C, E ) \ \to \ H^0 ( C, F ) \ \to \ H^0 ( C, F/E ) \ \to \\
 H^0 ( C, \Kc \otimes E^* )^* \ \to \ H^0 ( C, \Kc \otimes F^* )^* \ \to \ 0 , \end{multline*}
from which it follows that
\begin{equation} h^0 ( C, F ) = h^0 ( C , E ) + f \text{ if and only if } h^0 ( C, F^* \otimes \Kc ) = h^0 ( C, E^* \otimes \Kc ) - e + f . \label{differentf} \end{equation}
Thus, generalising (\ref{LineBundleHefQef}), we have a Cartesian diagram
\begin{equation} \xymatrix{ \HerS_\nd \ar[r]^\alpha & \Quot^{0, e} ( E^* ) \ar@{-->}[r]^a & \Urde \\
 \Hef ( \opeo )_\nd \ar[r]^-\alpha \ar@{^{(}->}[u] & \Qef \left( E , \Kc , H^0 ( C, E^* \otimes \Kc ) \right) \ar@{-->}[r]^-a \ar@{^{(}->}[u] & B_{r, d+e}^{h + f} \ar@{^{(}->}[u] } , \label{VectorBundleHefQef} \end{equation}
where $\alpha$ is as defined in Proposition \ref{alphaProperties} and $a \left( [ F^* \to E^* ] \right) = F$. The maps $a$ and $a \circ \alpha$ are generalisations of the Abel--Jacobi map.

Now for any $F$ of rank $r$, we write $H^0 ( C, \Hom ( F^* , E^* ) )_\inj$ for the open subset of generically injective maps. Notice that there is a natural action of $\Aut ( F^* )$ on $H^0 ( C, \Hom ( F^* , E^* ) )$ preserving $H^0 ( C, \Hom ( F^* , E^* ) )_\inj$.

\begin{proposition} \label{AJ} Let $E \to C$ be a bundle of rank $r$ and degree $d$.
\begin{enumerate}
\item[(a)] For a fixed bundle $F$ of rank $r$ and degree $d + e$, the locus
\[ \{ [ F_1^* \xrightarrow{j} E^* ] \in \Quot^{0, e} ( E^* ) : F_1 \cong F \hbox{ as vector bundles} \} \ =: \ Q_E (F) \]
is in bijection with $H^0 ( C, \Hom ( F^* , E^* ) )_\inj / \Aut ( F^* )$.
\item[(b)] Suppose that $\Quot^e ( E^* )$ contains a point $[ F^* \to E^* ]$ such that $F$ is a stable vector bundle. Let $a$ be as defined in (\ref{VectorBundleHefQef}). Then
\[ a^{-1} ( F ) \ = \ Q_E ( F ) \ \cong \ \PP H^0 ( C, \Hom ( F^* , E^* ))_\inj . \] \end{enumerate} \end{proposition}
\begin{proof} (a) The association $j \mapsto [ F^* \xrightarrow{j} E^* ]$ defines a map
\[ H^0 ( C, \Hom ( F^* , E^* ) )_\inj \ \to \ Q_E (F) \]
which is clearly surjective. 
 By definition of the Quot scheme and the universal property of kernels, two maps $j'$ and $j$ define the same element of $\Quot^{0, e} ( E^* )$ if and only if $j' = j \circ \gamma$ for some $\gamma \in \Aut ( F^* )$. This proves (a). Part (b) follows from (a) since $\Aut (F^*) = \C^*$ if $F$ is stable. \end{proof}

\begin{remark} If $L$ is a line bundle of degree $e$, then by Proposition \ref{AJ} we obtain
\[ a^{-1} ( L ) \ \cong \ Q_\Oc ( L ) \ \cong \ \PP H^0 ( C , \Hom ( L^{-1} , \Oc ) )_\inj \ \cong \ | L | . \]
Thus it makes sense to call $Q_E (F)$ a ``generalised linear series of $F$ with respect to $E$''. If $\Omega \subseteq H^0 ( C, \Hom ( F^* , E^* ) )$ is an $\Aut (F^*)$-invariant subspace, then the locus $\Omega_\inj / \Aut (F^*)$ generalises the notion of possibly incomplete linear subseries. \end{remark}

\begin{remark} 
The Riemann--Kempf singularity theorem \cite[Chap.\ VI.2]{ACGH}, gives a description of the tangent cones to the rank one Brill--Noether loci $W^f_e$. For $k \ge r$ and $W \in B_{r, d}^k$ a generated vector bundle, there is a generalisation of the Riemann--Kempf singularity theorem given as follows. For $\ell \ge r$, let $U^\ell \subseteq \Gr ( \ell , H^0 ( C, W ) )$ be the open subset of subspaces which generically generate $W$. Assuming this is nonempty, by \cite[Theorem 5.6]{Hit20} the projectivised tangent cone to $B_{r, d}^k$ at $W$ is given by
\begin{equation} \PP \cT_W B_{r, d}^k \ = \ \overline{ \bigcup_{\Lambda \in U^k} \left( \bigcap_{\Pi \in \left( U^r \cap \Gr ( r, \Lambda ) \right)} \Span \, \psi \left( Z_{\Pi} \times_C \PP W \right) \right) } . \label{RKSexpr} \end{equation}
Here $\psi \left( Z_\Pi \times_C \PP W \right)$ is a certain generalisation of the canonical image of the divisor associated to a section of a line bundle, defined in \cite[{\S} 5.2]{Hit20}. The point of interest for us here is that $U^r \cap \Gr ( r, \Lambda )$ appears in (\ref{RKSexpr}) in place of the linear subseries $g_d^{k-1}$ in the original Riemann--Kempf singularity theorem \cite[p.\ 241]{ACGH}. Continuing from the previous remark, we will interpret $U^r \cap \Gr ( r, \Lambda)$ as a generalised linear series.

Firstly, the map $( s_1 , \ldots , s_r ) \mapsto \Span \{ s_1 , \ldots , s_r \}$
defines a morphism
\[ H^0 ( C, \Hom ( \Oc^{\oplus r} , W ) )_\inj \ \to \ \Gr ( r, H^0 ( C, W ) ) , \]
whose image is exactly $U^r$, and which descends to an isomorphism
\begin{equation} H^0 ( C , \Hom ( \Oc^{\oplus r} , W ) )_\inj / \GL_r \ \isom \ U^r . \label{plusrIsomGr} \end{equation}
 By Proposition \ref{AJ} (a), this is exactly the generalised linear series $Q_{W^*} ( \Oc^{\oplus r} )$. Moreover, for $\Lambda \in U^k$ the isomorphism (\ref{plusrIsomGr}) restricts to an isomorphism
\[ \Lambda^{\oplus r}_\inj / \GL_r \ \isom \ U^r \cap \Gr ( r, \Lambda ) . \]
Thus $U^r \cap \Gr ( r, \Lambda )$ does indeed play the role of a linear series in the generalised Riemann--Kempf theorem. \end{remark}

\section{Some nonemptiness results} \label{nonempt}

In this section we give various sufficient or necessary conditions for nonemptiness of $\HefV_\nd$ and $\QefV$, for some values of $e$ and $f$. We begin with an easy ``theoretical bound'' valid for any smooth nondegenerate variety $S$ of dimension at least $2$ (it is trivial for curves).

\begin{lemma} \label{TheoreticalBound} Let $S$ be a smooth variety of dimension $r$ and $\ell = (\cL , V)$ a very ample linear series of dimension $n$ on $S$. Set $e_0 := n - r + 2$. Then for any $e > e_0$ and $f \le e - e_0$, the secant locus $\Hefl$ is nonempty. \end{lemma}

\begin{proof} By Theorem \ref{KleimanTheorem}, for a general $\Pi = \PP^{n - r + 1} \subset \PP V^*$ the intersection $\Pi \cap \psi (S)$ is smooth of dimension $1$. 
Any reduced subscheme $Z \subset \left( \Pi \cap \psi (S) \right)$ of length $e > e_0$ has linear span contained in $\Pi$, so has defect at least 
 $e - e_0$. Thus $Z$ defines a point of $\Hefl$ for any $f \le e - e_0$. \end{proof}

\begin{remark} If $e = n$, then $e - e_0 = r - 2$. Thus for $r \ge 3$, Lemma \ref{TheoreticalBound} implies that $H^{n-f}_n ( \ell )$ is nonempty for $1 \le f \le r-2$. This generalises the nonemptiness statement of \cite[Lemma 2.1]{AS} for $r \ge 3$. \end{remark}

Let us now focus on secant loci of scrolls over curves. Firstly, we have a technical statement.

\begin{lemma} \label{IdealSection} Let $\pi \colon \PP E \to C$ be a projective bundle, and $\iota \colon C \to \PP E$ a section. Let $x \in C$ be a point.
\begin{enumerate}
\item[(a)] There exists a local parameter $z$ on $C$ at $x$ and a frame $\phi_1 , \ldots , \phi_r$ for $E^*$ near $x$ such that for $k \ge 1$ the one-point subscheme $\iota ( kx ) \subset \PP E$ is defined by the ideal
\begin{equation} \left( \pi^* z^k , \frac{\phi_2}{\phi_1} , \ldots , \frac{\phi_r}{\phi_1} \right) . \label{KernelIdeal} \end{equation}
\item[(b)] Write $Z := \iota ( kx )$ and set $E_Z^* = \pi_* \cI_Z (1)$. Then $\{ \pi^* z^k \cdot \phi_1 , \phi_2 , \ldots , \phi_r \}$ is a frame for $E_Z^*$ near $x$.
\end{enumerate}
\end{lemma}

\begin{proof} Let $U := \Spec A$ be an affine neighbourhood of $x$ in $C$ upon which $E$ is trivial. Let $\nu_1 , \ldots , \nu_r$ be a frame such that $\iota \colon C \to \PP E$ corresponds to the subbundle $\Oc \cdot \nu_1$. Let $\phi_1 , \ldots , \phi_r$ be the dual frame of $E^*|_U$, and let $z$ be a local parameter on $C$ at $x$. We consider the affine variety
\[ \Omega \ := \ \Spec A \left[ \frac{\phi_2}{\phi_1} , \ldots , \frac{\phi_r}{\phi_1} \right] \ \cong \ U \times \bA^{r-1} \ \subset \ \PP E . \]
In these coordinates, the embedding
\[ kx \ = \ \Spec \left( A / \fm_x^k \right) \ \hookrightarrow \ U \ \xrightarrow{\iota} \ \Omega \]
corresponds to the composed map of rings
\[ A \left[ \frac{\phi_2}{\phi_1} , \ldots , \frac{\phi_r}{\phi_1} \right] \ \xrightarrow{\sigma^*} \ A \ \to \ A / \fm_x^k . \]
where $\sigma^* z = z$ and $\sigma^* ( \phi_i / \phi_1 ) = 0$ for $2 \le i \le r$. The last composed map has kernel exactly (\ref{KernelIdeal}), and we obtain (a).

For the rest: Using (a), we see that $\cI_{\iota ( ex )}(1)$ is generated on $\pi^{-1} ( U )$ by the elements $\pi^* z^k \cdot \phi_1 , \phi_2 , \ldots , \phi_r$. Statement (b) follows. 
\end{proof}

\begin{corollary} Let $\PP E$ and $\iota$ be as above. For any effective divisor $D$ of degree $e$ on $C$, the subscheme $\iota (D) \subset \PP E$ defines a point of $\Hilb^e ( \PP E )$ which is $\pi$-nondefective and smoothable. \label{NondefSection} \end{corollary}

\begin{proof} Firstly, as $\iota (D)$ is the image of a curvilinear scheme by an embedding, $\iota (D)$ is smoothable. For the rest: As $\pi$-nondefectivity is local on $C$, it suffices to consider the case $D = ex$. In this case, by Lemma \ref{IdealSection} (b) we see that 
\[ E_Z^*|_U \ = \ \pi_* \cI_Z ( 1 )|_U \ \cong \ \cO_U ( -e x ) \oplus \cO_U^{\oplus (r-1)} . \]
Therefore, $E^* / E_Z^*$ has length $e$; whence $\iota (D)$ is $\pi$-nondefective by Definition \ref{DefnPiDef}. \end{proof}

Suppose now that $S = \PP E$ is a scroll of dimension $r \ge 2$ over $C$. We will use existing results on secant loci of curves to prove nonemptiness of secant loci of $S$ whose expected dimension is large.

\begin{proposition} Let $E$ be a bundle of rank $r \ge 2$ and degree $d$, and write $S := \PP E$. Suppose that $V \subseteq H^0 ( S, \cL_M )$ is a linear subsystem of dimension $n + 1$. 
 Suppose that
\begin{equation} 0 \ < \ f \ < \ e \quad \hbox{and} \quad re - f ( n + 1 - e + f ) \ \ge \ (r-1)e . \label{dimhyp} \end{equation}
Then $\HefV_\nd$ and $\QefV$ are nonempty. \label{nonempty} \end{proposition}

\begin{proof} 
Write $d' := \deg (E \otimes M^{-1} ) = d - r \cdot \deg (M)$. By \cite[Proposition 6.1]{PR}, for all $k \gg 0$, there exists a short exact sequence $0 \to L \to E \otimes M^{-1} \to F \to 0$ where $F$ is a stable vector bundle of rank $r-1$ and degree $k$. We choose
\begin{equation} k \ \ge \ \max\{ 0, 2g + d' \} . \label{soineq} \end{equation}
Since $k \ge 0$, we have $h^0 ( C, F^* ) = 0$ by stability. Taking global sections of $0 \to F^* \to E^* \otimes M \to L^{-1} \to 0$ and dualising, we obtain
\begin{equation} \cdots \ \to \ H^0 ( L^{-1} )^* \ \to \ H^0 ( E^* \otimes M )^* \ \to \ 0 . \label{LtoEMsurj} \end{equation}
Projectivising and identifying $C$ with $\PP L$, we obtain a commutative diagram
\[ \xymatrix{ | L^{-1} |^* \ar@{-->>}[r] & \PP V^* \\
C \ar[u]_{\varphi_{L^{-1}}} \ar@{^{(}->}[r]^\iota & \PP E . \ar@{-->}[u]} \]
Now $\deg L^{-1} = 
 k - d'$, which by (\ref{soineq}) is at least $2g$. 
 Hence $h^1 ( C, L^{-1} ) = 0$ and $\varphi_{L^{-1}}$ is an embedding. As $C$ is nondegenerate in $| L^{-1} |^*$, we obtain an identification of $\PP V$ with a subsystem $g_{k - d'}^n$ of $| L^{-1} |$.

We recall next from \cite[p.\ 355]{ACGH} that the \textsl{relative deficiency} of this $g_{k - d'}^n$ is defined to be $\deg L^{-1} - g - n$. As $h^1 ( C, L^{-1} ) = 0$, we have
\[ \deg L^{-1} - g - n \ 
 = \ \chi ( C, L^{-1} ) - ( n + 1 ) \ = \ h^0 (C, L^{-1} ) - \dim V , \]
which is nonnegative since we have identified $\PP V$ with a subsystem of $| L^{-1} |$ in (\ref{LtoEMsurj}). Furthermore, by (\ref{dimhyp}) we have $e - f(n + 1 - e + f) \ge 0$. Hence by loc.\ cit.\ the usual secant locus $\Vef ( g_{d + k}^n ) \subseteq \Ce$ is nonempty. Thus there exists a divisor $D$ on $C$ of degree $e$ such that $\Span \, \psi ( \iota (D) )$ has dimension $e-f-1$. Then $\iota (D)$ is a length $e$ subscheme of $\PP E$ spanning at most a $\PP^{e-f-1}$ in $|\cL_M|^*$; and $\iota(D)$ is $\pi$-nondefective and smoothable by Corollary \ref{NondefSection}. Hence $\HefV_\nd$ is nonempty, and so therefore is $\QefV$ by Proposition \ref{HQequiv}. \end{proof}

\begin{remark} If $f = 1$, then (\ref{dimhyp}) reduces to $e \ge 2$ and $2e - 1 \ge n + 1$. 
If $e = n$, then this is satisfied for $n \ge 2$. 
 Thus for $r \ge 2$, the statement of Proposition \ref{nonempty} again generalises the nonemptiness part of \cite[Lemma 2.1]{AS}.
\end{remark}

We conclude this section with a necessary condition for nonemptiness.

\begin{proposition} Let $E$ be a bundle of rank $r$ and degree $d$. Set $\hh := h^0 ( C, \Kc \otimes E )$. If $\Qef$ is nonempty, then the twisted Brill--Noether locus
\[ B^{\hh + f}_{1, e} ( \Kc \otimes E ) \ = \ \{ L \in \Pic^e ( C ) : h^0 ( C, \Kc L \otimes E ) \ge f \} \]
contains an effective line bundle. \end{proposition}

\begin{proof} Suppose $[ F^* \to E^* ]$ belongs to $\Qef$. Using (\ref{differentf}), we obtain
\[ h^0 ( C, \Kc \otimes F ) \ \ge \ h^0 ( \Kc \otimes E ) + f . \]
 Now $F \subseteq E (D)$ for some $D \in \Ce$. Then $\Oc (D)$ belongs to $B^{\hh + f}_{1, e} ( \Kc \otimes E )$. \end{proof}

\noindent This observation will be further developed in {\S} \ref{param} and thereafter.

\section{Tangent spaces to secant loci} \label{TangSp}

In this section, we will describe the Zariski tangent spaces of $\Qef (V)$, generalising \cite[Lemma IV.1.5]{ACGH} and \cite[Theorem 0.3]{Cop}. In preparation, we give a construction of the first-order infinitesimal deformations of an element of any Quot scheme, which will be convenient for calculations.

\subsection{Tangent spaces to Quot schemes}

Let $( X, \Ox (1) )$ be a polarised variety and $\cV$ a coherent sheaf over $X$. The scheme $\Quot^\Phi ( \cV )$ parametrises equivalence classes of coherent quotients $\cV \xrightarrow{q} \cQ$ where $\cQ$ has Hilbert polynomial $\Phi$. The first-order infinitesimal deformations of a given $[ \cV \xrightarrow{q} \cQ]$ are parametrised by the Zariski tangent space $T_q \Quot^\Phi ( \cV )$. This is canonically isomorphic to $H^0 ( X, \Hom ( \cW , \cQ ))$, where $\cW = \Ker (q)$; see for example \cite[Chapter 2]{HL}. These deformations naturally induce deformations of the sheaf $\cW$, which we will now construct. (A similar approach to deformations is given in \cite[Exercise 6.12 c]{Eis}.)

For each $u \in H^0 ( X, \Hom ( \cW, \cQ ))$, let $\bW_u$ be the sheaf over $\Spec \C [\varepsilon] \times X$ given over each open set $U \subseteq X$ by
\begin{equation} \bW_u ( U ) \ = \ \{ \varepsilon t + s \in H^0 ( U, \varepsilon \cdot \cV \oplus \cW ) : q(t) = u(s) \} \label{defFu} \end{equation}
(Here we use the fact that $\Spec \, \C[\varepsilon] \times X$ has the same ambient topological space as $X$.) The sheaf $\bW_u$ is naturally contained in $\varepsilon \cdot \cV \oplus \cV$, and $\bW_u = \bW$ is easily checked to fit into an exact diagram of the form
\begin{equation} \xymatrix{ & 0 \ar[d] & 0 \ar[d] & 0 \ar[d] & \\
0 \ar[r] & \varepsilon \cdot \cW \ar[d]^j \ar[r] & \bW \ar[r] \ar[d] & \cW \ar[r] \ar[d]^j & 0 \\
 0 \ar[r] & \varepsilon \cdot \cV \ar[d]^q \ar[r] & \varepsilon \cdot \cV \oplus \cV \ar[r] \ar[d]^{q'} & \cV \ar[r] \ar[d]^q & 0 \\
 0 \ar[r] & \varepsilon \cdot \cQ \ar[r] \ar[d] & \bQ \ar[r] \ar[d] & \cQ \ar[r] \ar[d] & 0 \\
& 0 & 0 & 0. & } \label{TangSpDiag} \end{equation}
Thus $\bW_u$ defines a first-order infinitesimal deformation of $\cV \xrightarrow{q} \cQ$.

\begin{lemma} \label{ExplDef} The association
\[ u \ \mapsto \ \left[ \varepsilon \cdot \cV \oplus \cV \to \frac{\varepsilon \cdot \cV \oplus \cV}{\bW_u} \right] \]
defines a bijection between $H^0 ( C, \Hom ( \cW, \cQ ) )$ and the set of isomorphism classes of first-order infinitesimal deformations of $[\cV \xrightarrow{q} \cQ]$. \end{lemma}

\begin{proof} We indicate only the main steps. Recall firstly that the \textsl{graph} of a map of $\Ox$-modules $u \colon \cW \to \cQ$ is the $\Ox$-submodule $\Gamma_u \subseteq \cQ \oplus \cW$ given over open $U \subseteq X$ by
\[ \Gamma_u ( U ) \ = \ \{ ( u(s), s ) : s \in \cW ( U ) \} . \]
By construction, $\bW_u$ is exactly the inverse image of $\Gamma_{\varepsilon u}$ by the projection $\varepsilon \cdot \cV \oplus \cW \to \varepsilon \cdot \cQ \oplus \cW$.

Now by for example \cite[{\S} 2.2]{HL}, each first-order infinitesimal deformation of $\cV \xrightarrow{q} \cQ$ corresponds to an exact diagram of the form (\ref{TangSpDiag}).
 As mentioned above, $\bW = \bW_u$ fits into such a diagram. Conversely, given a diagram of the form (\ref{TangSpDiag}), write $\tbW := \left( q \oplus \Iden_\cV \right) ( \bW )$. Then
\begin{enumerate}
\item[(i)] $\tbW \cap ( \varepsilon \cdot \cQ \oplus 0 ) = 0$ since $\bW \cap ( \varepsilon \cdot \cV \oplus 0 ) = \varepsilon \cdot \cW$; and
\item[(ii)] the image of $\tbW$ in $\cV$ is $\cW$ since $\bW \to \cW$ is surjective.
\end{enumerate}
Linear algebra arguments then show that $\tbW = \Gamma_{\varepsilon u}$ for a uniquely determined $u \colon \cW \to \cQ$, and that that $\bW = \bW_u$. \end{proof}

\subsection{Zariski tangent spaces of \texorpdfstring{$\QefV$}{Q e-f e (V)}}

We now consider a smooth curve $C$ and a vector bundle $E \to C$ of rank $r$ and degree $d$. Let $M$ be a line bundle and $V \subseteq H^0 ( C, E^* \otimes M )$ a subspace of dimension $n + 1$. We will study the tangent spaces to $\Qef ( E, M, V ) =: \QefV$.

For the remainder of the paper, for any sheaf $\cW$ over $C$ we abbreviate $H^i ( C, \cW )$ and $h^i ( C, \cW )$ to $H^i ( \cW )$ and $h^i ( \cW )$ respectively.

Following the treatment in \cite[Chap.\ IV]{ACGH} of the tangent spaces of $\Ce^f$, which is exactly $\Qef ( \Kc , \Oc , H^0 ( \Kc ) )$, we fix the following. Consider an elementary transformation $0 \to E \to F \to \cT \to 0$ where $\tau$ has length $e$. Taking $\Hom ( - , M )$, we obtain a sequence
\begin{equation} 0 \ \to \ F^* \otimes M \ \to \ E^* \otimes M \ \xrightarrow{q} \ \tau \ \to \ 0 \label{hEFt} \end{equation}
where $\tau := Ext^1_\Oc ( \cT , M )$ is noncanonically isomorphic to $\cT$. By the previous section, and since $M$ is invertible,
\[ T_{F^*} \Quot^{0, e} ( E^* ) \ \cong \ H^0 ( C, \Hom ( F^*, \tau \otimes M^{-1} ) ) \ \cong \ H^0 ( C, \Hom ( F^* \otimes M , \tau ) ) . \]
Also, recall that there is a natural map
\[ c \colon H^0 ( \Hom ( F^* \otimes M , \tau ) ) \ \to \ \Hom \left( H^0 ( F^* \otimes M ) , H^0 ( \tau ) \right) . \]

\begin{proposition} \label{TangSpGen} Let $E$, $M$ and $V$ be as above. Suppose that $[ F^* \to E^* ]$ is a point of $\QefV \setminus Q^{e-f-1}_e ( V )$. Then
\[ T_{F^*} \QefV \ = \ \{ u \in H^0 ( \Hom ( F^* \otimes M , \tau )) : c(u) ( V \cap H^0 ( F^* \otimes M )) \subseteq q ( V ) \} . \]
Moreover, this is precisely the kernel of the map
\[ H^0 ( \Hom ( F^* \otimes M , \tau ) ) \ \to \ \Hom \left( V \cap H^0 ( F^* \otimes M ) , \frac{H^0 ( \tau )}{q(V)} \right) \]
given by composing $c$, restriction to $V \cap H^0 ( F^* \otimes M )$ and quotient by $q(V)$.
\end{proposition}

\begin{proof} A tangent vector $u \in H^0 ( \Hom ( F^* \otimes M , \tau ))$ belongs to $T_{F^*} \QefV$ if and only if every $s \in V \cap H^0 ( F^* \otimes M )$ lifts to the deformation $\bW_u$ given as in (\ref{defFu}) by
\[ \bW_u ( U ) \ = \ \left\{ \varepsilon t + s \in \varepsilon \cdot ( E^* \otimes M ) (U) \oplus ( F^* \otimes M ) (U) : q ( t ) = u ( s ) \right\} \]
over each open subset $U \subseteq C$. This is equivalent to saying that for each element $s \in V \cap H^0 ( F^* \otimes M )$, there exists $t \in V$ such that $q (t) = c(u) (s)$; in other words, that $c(u) \left( V \cap H^0 ( \hF ) \right) \subseteq q ( V )$. This proves the first statement; and the second is then clear. \end{proof}

\subsection{The complete case}

During this subsection, we fix $M = \Kc$ and assume $V = H^0 ( E^* \otimes \Kc )$. In this case, we can say more about $T_{F^*} \Qef$.

Firstly, we recall some familiar facts. To ease notation, for any locally free sheaf $W$ over $C$ we will denote $\Kc \otimes W^*$ by $\widehat{W}$. Fix an elementary transformation $0 \to E 
 \to F \to \cT \to 0$. This canonically determines a sequence $0 \to \hF \to \hE \xrightarrow{q} \tau \to 0$ and an element of $\Quot^{0, e} ( E^* )$ as above. Let $\partial \colon H^0 ( \tau ) \to H^1 ( \hF )$ denote the associated coboundary map.  Then we obtain a commutative diagram with exact rows
\begin{equation} \xymatrix{ H^0 ( \Hom ( \hF , \tau ) ) \ar[d]^c \ar[r]^\delta & H^1 ( \End \hF ) \ar[d]^\cup \ar[r]
 & H^1 ( \Hom ( \hF , \hE ) ) \ar[d] \ar[r] & 0 \\
 \Hom \left( H^0 ( \hF ) , H^0 ( \tau ) \right) \ar[r]^{\partial_*} & \Hom \left( H^0 ( \hF ) , H^1 ( \hF ) \right) \ar[r] & \Hom \left( H^0 ( \hF ) , H^1 ( \hE ) \right) \ar[r] & 0 } \label{CompleteTangSp} \end{equation}
where $\cup$ is the cup product map. 
 We recall the well known fact that via Serre duality, $\cup$ is dual to the Petri map
\[ \mu \colon H^0 ( \hF ) \otimes H^0 ( F ) \ \to \ H^0 ( \hF \otimes F ) \ = \ H^0 ( \Kc \otimes \End F ) . \]

\begin{proposition} \label{CompleteZarTangSp} Let $0 \to E \to F \to \cT \to 0$ be as above. Suppose that $[ F^* \to E^* ]$ belongs to $\Qef \setminus Q^{e-f-1}_e$. 
\begin{enumerate}
\item[(a)] We have $T_{F^*} \Qef = \Ker ( \cup \circ \delta ) \ = \ \Image \left( ^t\delta \circ \mu \right)^\perp$.
\item[(b)] The locus $\Qef$ is smooth and of the expected dimension at $[ F^* \to E^* ]$ if and only if $\mu^{-1} \left( H^0 ( \hF \otimes E ) \right) = H^0 ( \hF ) \otimes H^0 ( E )$.
\end{enumerate}
\end{proposition}

\begin{proof} (a) Since $V = H^0 ( \hE )$, we have $\Image ( c(u) ) \subseteq q ( V )$ if and only if $\partial \circ c (u)$ is zero in $\Hom \left( H^0 ( \hF ) , H^1 ( \hF ) \right)$.
 Hence by Proposition \ref{TangSpGen} and commutativity of (\ref{CompleteTangSp}), it follows that
\[ T_{F^*} \Qef \ = \ \Ker ( \partial_* \circ c ) \ = \ \Ker ( \cup \circ \delta ) . \]
As $\Ker ( \varphi ) = \Image ( {^t\varphi} )^\perp$ for any vector space map $\varphi$, and since ${^t\cup} = \mu$, this in turn coincides with $\Image \left( ^t\delta \circ \mu \right)^\perp$.

(b) By part (a), the codimension of $T_{F^*} \Qef$ in $T_{F*} \Quot^e ( E^* )$ is exactly $\dim \Image ( {^t\delta} \circ \mu )$.
 Let us compute the latter. Dualising (\ref{CompleteTangSp}), we obtain
\begin{equation} \label{DualCompleteTangSp} \xymatrix{ 0 \ar[r] & H^0 ( \hF ) \otimes H^0 ( E ) \ar[d] \ar[r] & H^0 ( \hF ) \otimes H^0 ( F ) \ar[r] \ar[d]_\mu & H^0 ( \hF ) \otimes H^0 ( \tau )^* \ar[d] \\
 0 \ar[r] & H^0 ( \hF \otimes E ) \ar[r] & H^0 ( \hF \otimes F ) \ar[r]^-{^t\delta} & H^0 ( \Hom ( \hF , \htau ) )^* . } \end{equation}
Thus $\Ker ( {^t\delta} \circ \mu ) = \mu^{-1} \left( H^0 ( \hF \otimes E ) \right)$, and so $\Image ( {^t\delta} \circ \mu )$ has dimension
\begin{multline*} h^0 ( \hF ) \cdot h^0 ( F ) - \dim \mu^{-1} \left( H^0 ( \hF \otimes E ) \right) \ = \\ 
 h^0 ( \hF ) \cdot \left( h^0 ( F ) - h^0 ( E )  \right) - \dim \frac{\mu^{-1} \left( H^0 ( \hF \otimes E ) \right)}{H^0 ( \hF ) \otimes H^0 ( E )} . \end{multline*}
By hypothesis, $h^0 ( \hF ) = h^0 ( \hE ) - e + f$. Hence by (\ref{differentf}) we have $h^0 ( F ) - h^0 ( E ) = f$. 
 Thus $\Image ( {^t\delta} \circ \mu )$ has dimension
\[ \left( h^0 ( \hE ) - e + f \right) \cdot f - \dim \frac{\mu^{-1} \left( H^0 ( \hF \otimes E ) \right)}{H^0 ( \hF ) \otimes H^0 ( E )} . \]
As the expected codimension of $\Qef$ in $\Quot^{0, e} ( E^* )$ is $( h^0 ( \hE ) - e + f ) \cdot f$ (for in this case $n + 1 = h^0 ( \hE )$), we see that $T_{F^*} \Qef$ has the expected codimension if and only if $\mu^{-1} \left( H^0 ( \hF \otimes E ) \right) = H^0 ( \hF ) \otimes H^0 ( E )$. \end{proof}

\subsection{Some examples} \label{examples}

Continuing to use the notation of the previous subsection, let us return to the study of the generalised Abel--Jacobi map
\[ a \colon \Quot^e ( E^* ) \ \dashrightarrow \ U_C ( r, d + e ) \]
given by $[ F^* \to E^* ] \mapsto F$. As noted in {\S} \ref{HefQefBN}, we have 
\[ a^{-1} \left( B_{r, d + e}^{h^0 ( E ) + f} \right) \ = \ \Qef ( E , \Kc , H^0 ( E^* \otimes \Kc ) ) . \]
It is well known that $B_{r, d + e}^{h^0 ( E ) + f}$ is smooth and of the expected dimension at a stable bundle $F$ if and only if $\mu$ is injective. In the special case $E = \Oc$ and $F = \Oc (D)$ with $D \in \Ce^f \setminus \Ce^{f+1}$, by \cite[Lemma IV.1.6]{ACGH}, injectivity of $\mu$ is also equivalent to smoothness of $\Ce^f$ at $D$. One can check that in this situation, the condition in Proposition \ref{CompleteZarTangSp} (b) is equivalent to injectivity of $\mu$, so there is no contradiction. 
 However, we will now give examples showing that in general, smoothness of $B_{r, d+e}^{h^0 ( E ) + f}$ at $F$ is neither necessary nor sufficient for smoothness of $\Qef ( E, \Kc , H^0 ( E^* \otimes \Kc ) )$ at $F^*$. 

\begin{example} \label{ExaOne} Let $D$ be an effective divisor of degree $e \le g - 1$, and $E := \Oc ( -D )$. Set $F = E (D) = \Oc$, so that $\hF = \Kc$. Here $\deg(E) = - e$, and $f = 1$, and
\[ F \ \in \ \Qeo ( E, \Kc , H^0 ( E^* \otimes \Kc ) ) \ = \ \Qeo ( \Oc(-D) , \Kc , H^0 ( \Kc ( D ) ) ) . \]
Trivially, $B_{1 , d+e}^1 = B_{1 , 0}^1 = \{ \Oc \}$ is smooth and of the expected dimension at $F$.

Now $H^0 ( \hF ) \otimes H^0 ( F ) = H^0 ( \Kc ) \otimes H^0 ( \Oc )$. Clearly $\mu$ is an isomorphism, and
\[ \mu^{-1} \left( H^0 ( \hF \otimes E ) \right) \ = \ \mu^{-1} \left( H^0 ( \Kc (-D) ) \right) \ = \ H^0 ( \Oc ) \otimes H^0 ( \Kc (-D) ) , \]
which is nonzero since $\deg (D) \le g-1$. On the other hand, $H^0 ( \Oc (-D) ) \otimes H^0 ( \Kc )$ is zero. By Proposition \ref{CompleteZarTangSp} (b), the secant locus $\Qeo$ is not smooth at $F$. (One can check that the expected dimension of $\Qeo$ is $e - g < 0$.)
\end{example}

The next example shows that $\Qef$ can be smooth at $F^*$ even though $B_{r, d}^k$ is not smooth at $F$. The idea is to produce an $E$ and $F$ such that $\Ker \, \mu$ is nonzero but contained in $H^0 ( \hF ) \otimes H^0 ( E )$ so that the condition $\mu^{-1} H^0 ( \hF \otimes E ) = H^0 ( \hF ) \otimes H^0 ( E )$ can still obtain. Perhaps somewhat artificially, we will start with the bundle $F$ and choose an appropriate $E$. The bundle $F$ will have rank two and canonical determinant; the Brill--Noether loci of such bundles have been studied in \cite{TiB04}, \cite{TiB08} and elsewhere, and are known to have excess dimension in many cases.

\begin{example} \label{ExaTwo} Suppose $g \ge 6$, so that $B_{1, g-2}^2$ is nonempty. 
 Let $L_1$ be a line bundle of degree $g-2$ with $h^0 ( L_1 ) = 2$. Let $x_1$ and $x_2$ be general points of $C$, and set $L_2 := \Kc L_1^{-1} ( - x_1 - x_2 )$. Then $h^0 ( L_2 ) = 1$. 
 Let
\begin{equation} 0 \ \to \ L_1 \oplus L_2 \ \to \ F \ \to \ \cO_{x_1} \oplus \cO_{x_2} \ \to \ 0 \label{FfirstET} \end{equation}
be an elementary transformation such that neither $L_1$ nor $L_2$ is a quotient of $F$. By \cite[Th\'eor\`eme A.5]{Mer}, we may assume that $F$ is stable.

We now claim that $H^0 ( F ) = H^0 ( L_1 ) \oplus H^0 ( L_2 )$. It will suffice to show that the composed map
\begin{equation} H^0 ( \cO_{x_1} \oplus \cO_{x_2} ) \ \to \ H^1 ( L_1 ) \oplus H^1 ( L_2 ) \ \to \ H^1 ( L_1 ) \ \cong \ H^0 ( \Kc L_1^{-1} )^* \label{AnotherComposedMap} \end{equation}
is injective. It is not hard to see that the image of this composed map is exactly the cone over the secant spanned by the points $x_1$ and $x_2$ on the plane curve $\varphi_{\Kc L_1^{-1}} (C)$. 
 As already used, by generality of $x_1$ and $x_2$ we have $h^0 ( \Kc L_1^{-1} ( - x_1 - x_2 ) ) ) = 1$, so this cone has dimension 2. Hence (\ref{AnotherComposedMap}) is an isomorphism, and in particular injective as desired.

Thus $F$ defines a point of $B^3_{2, 2g-2}$. Moreover, since $\det (F) \cong \Kc$, we have $F \cong \Kc \otimes F^*$, and the Petri map can be identified with the natural map
\[ H^0 ( F ) \otimes H^0 ( F ) \ \to \ H^0 ( F \otimes F ) . \]
Abusing notation, we also denote this map by $\mu$. Let $s_1$, $s_2$ be a basis for $H^0 ( L_1 )$. Then $s_1 \wedge s_2$ is a nonzero element of $\Ker ( \mu )$, so $B^3_{2, 2g-2}$ 
 fails to be smooth and/or of expected dimension at $F$.

We will now construct a bundle $E$ such that $[F^* \to E^*]$ is a smooth point of a certain $\Qef ( E^* , \Kc , H^0 ( E^* \otimes \Kc ) )$. Let $E$ be an elementary transformation $0 \to E \to F \to \cO_D \to 0$ where $D$ is a general effective divisor of degree $e \ge 2$, and
\begin{equation} \Image \left( E|_p \to F|_p \right) \ = \ L_1|_p \hbox{ for each } p \in D . \label{DefnE} \end{equation}
As $L_1 \to F \to \cO_D$ is zero, $L_1 \subset E$. (Note that $E$ is therefore not stable, as $\deg (L_1) = g - 2 \ge \frac{2g-2-e}{2}$, but this does not affect the argument.) In particular $H^0 ( E )$ contains $H^0 ( L_1 )$. Furthermore, let $s_3$ be a generator of $H^0 ( L_2 )$. By generality of $D$, we may assume that $s_3 (p) \not\in L_1|_p$ for 
 all $p \in D$. Thus $s_3$ is not a section of $E$, whence $h^0 ( E ) = h^0 ( L_1 ) = 2$. By Riemann--Roch, $h^0 ( E^* \otimes \Kc ) = e + 2$. 
 As
\[ h^0 ( F^* \otimes \Kc ) \ = \ h^0 ( F ) \ = \ 3 \ = \ 
h^0 ( E^* \otimes \Kc ) - e + 1 , \]
we have $[ F^* \to E^* ] \in \Qeo \left( E , \Kc , H^0 ( C, E^* \otimes \Kc ) \right)$. We will show that it is a smooth point.

By Proposition \ref{CompleteZarTangSp} (b), noting that $F \cong F^* \otimes \Kc = \hF$, we must show that
\begin{equation} \mu^{-1} \left( H^0 ( F \otimes E ) \right) \ = \ H^0 ( F ) \otimes H^0 ( E ) . \label{SmoothCond} \end{equation}

Let $s := \sum_{i, j=1}^3 \alpha_{ij} s_i \otimes s_j$ be an element of $H^0 ( F ) \otimes H^0 ( F )$, and assume that $( \mu(s) )(p)$ belongs to the image of $F \otimes E|_p$ for all $p \in C$. This is a nonempty condition only at $p \in D$, since $E|_p = F|_p$ for $p \not\in D$. For $p \in D$, using (\ref{FfirstET}), there is a canonical splitting of $(F \otimes F)|_p$ as
\[ \left( L_1|_p \otimes L_1|_p \right) \oplus \left( L_2|_p \otimes L_1|_p \right) \oplus \left( L_1|_p \otimes L_2|_p \right) \oplus \left( L_2|_p \otimes L_2|_p \right) . \]
We write out $\left( \mu(s) \right) ( p )$ in terms of this splitting:
\begin{multline} \left( \alpha_{11} s_1 (p) \otimes s_1 (p) + \alpha_{21} s_2 (p) \otimes s_1 (p) + \alpha_{12} s_1 (p) \otimes s_2 (p) + \alpha_{22} s_2 (p) \otimes s_2 (p) , \right. \\
 \alpha_{31} s_3 (p) \otimes s_1 (p) + \alpha_{32} s_3 (p) \otimes s_2 (p) , \alpha_{13} s_1 (p) \otimes s_3 (p) + \alpha_{23} s_2 (p) \otimes s_3 (p) , \\
\left. \alpha_{33} s_3 (p) \otimes s_3 (p) \right) . \label{sofp} \end{multline}
Now in view of (\ref{DefnE}), by hypothesis $\left( \mu(s) \right) (p)$ belongs to
\[ \Image \left( E|_p \otimes F|_p \ \to \ F|_p \otimes F|_p \right) \ = \ \left( L_1|_p \oplus L_2|_p \right) \otimes L_1|_p\]
for $p \in D$. Thus the third and fourth components of (\ref{sofp}) must be zero for all $p \in D$. Firstly, since $s_3 (p)$ can be assumed to be nonzero for all $p \in D$ by generality of $D$, we obtain $\alpha_{33} = 0$. Furthermore, considering the third component of (\ref{sofp}), we obtain
\[ \left( \alpha_{13} s_1 (p) + \alpha_{23} s_2 (p) \right) \otimes s_3 (p) \ = \ 0 \hbox{ for all } p \in D . \]
Using the nonvanishing of $s_3 (p)$ above, we obtain $\alpha_{13} s_1 (p) + \alpha_{23} s_2 (p) = 0$ for all $p \in D$. Again using generality of $D$, we may assume that 
 no two points of $D$ are identified under the map $C \to |L_1|^* = \PP^1$. Thus this condition places $e$ independent linear conditions on $(\alpha_{13}, \alpha_{23} )$. Since $e \ge 2$, the only solution is $\alpha_{13} = \alpha_{23} = 0$. It follows that $s \in H^0 ( F ) \otimes H^0 ( E )$. Thus (\ref{SmoothCond}) is satisfied and $\Qeo$ is smooth at $[ F^* \to E^* ]$. \end{example}

\begin{remark} In example 
 \ref{ExaOne}, the singularity arises from a determinantal variety which does not have the expected dimension. Moreover, the bundle $E$ in Example \ref{ExaTwo} is not semistable. It would be interesting to have examples of other kinds, and to investigate more thoroughly the relationships between smoothness of $B_{r, d}^k$, smoothness of $\Qef$ and semistability of $E$. \end{remark}

\section{Parametrisation of \texorpdfstring{$\Qeo$}{Q e-1 e}} \label{param}

As before, let $E$ be a vector bundle of rank $r$ and degree $d$ over $C$. In \cite[{\S} 4.3]{Hit19}, a parametrisation is given for the inflectional loci of a map $\PP E \dashrightarrow |\cL_M|^*$. We will now construct a similar parametrisation of $\Qeo ( E, M, H^0 ( E^* \otimes M ) )$. In the sections that follow, we will give some applications similar to those in \cite[{\S} 5--6]{Hit19} for this parametrisation.

For $M \in \Pic^0 (C)$, consider the Quot scheme
\[ \Quot^{r-1, d + 2r(g-1) + e} ( \Kc M^{-1} \otimes E ) \]
parametrising equivalence classes of quotients $[ \Kc M^{-1} \otimes E \to \cQ ]$ where $\cQ$ is coherent of rank $r-1$ and degree $\deg (\Kc M^{-1} \otimes E ) + e$; equivalently, invertible subsheaves of $\Kc M^{-1} \otimes E$ of degree $-e$. As we wish to focus on subsheaves, and to ease notation, we denote this Quot scheme by $\qoe$.

Let $a \colon \Ce \to \Pic^{-e} (C)$ be the Abel--Jacobi-type map $D \mapsto \Oc (-D)$, and let $b \colon \qoe \to \Pic^{-e} (C)$ be the forgetful map $\left[ L \xrightarrow{\sigma} \Kc M^{-1} \otimes E \right] \mapsto L$. Consider the fibre product
\begin{equation} \xymatrix{ S^e_M \ar[r] \ar[d] & \qoe \ar[d]^b \\
 \Ce \ar[r]^-a & \Pic^{-e} (C). } \label{SeMdefn} \end{equation}
Set-theoretically, $S^e_M$ parametrises pairs $( \sigma , D )$ where $D$ is an effective divisor of degree $e$ on $C$ and $\sigma \colon \Oc (-D) \to \Kc M^{-1} \otimes E$ is a sheaf injection. We will now construct a family of elementary transformations of $E^*$ parametrised by $S^e_M$.

Consider the following diagram where all maps are projections.
\[ \xymatrix{ & \qoe \times \Ce \times C \ar[dl]^{p_1} \ar[dd]_{q_0} \ar[dr]_{p_2} & \\
 \qoe \times C \ar[dr]_{q_1} & & \Ce \times C \ar[dl]^{q_2} \\
 & C & } \]
Recall also the exact sequence $0 \to \cO_{\Ce \times C} ( - \scrD ) \to \cO_{\Ce \times C} \xrightarrow{\rho} \cO_\scrD \to 0$ where $\scrD$ is the universal divisor.

Let $\Sigma \colon \cL \to q_1^* ( \Kc M^{-1} \otimes E )$ be the universal subsheaf over $\qoe \times C$. A key point is that $\cL$ is locally free, and hence $p_1^* \cL$ is too (this is the reason we work with $\qoe$ instead of $\Quot^{1, e} ( \Kc^{-1} M \otimes E^* )$ directly). Thus $p_1^* \Sigma$ can be considered as a map of vector bundles over $\qoe \times \Ce \times C$ (with a degeneracy locus). Hence we can take its transpose, twist by $q_0^* ( \Kc M^{-1} )$ and form the composed map
\begin{equation} q_0^* E^* \ \to \ q_0^* ( \Kc M^{-1} ) \otimes p_1^* \cL^\vee \ \to \ q_0^* ( \Kc M^{-1} ) \otimes p_1^* \cL^\vee \otimes p_2^* \cO_{\scrD} . \label{StillAnothCompMap} \end{equation}
Now $\Ce$ is also a Quot scheme, parametrising degree $e$ torsion quotients of any line bundle over $C$. Thus, by the universal property of Quot schemes, the restriction of this map to any fibre
\[ \left\{ [ L \xrightarrow{\sigma} \Kc M^{-1} \otimes E ] , D \right\} \times C \]
is identified with the composed map
\[ E^* \ \xrightarrow{\hsigma} \ \Kc M^{-1} L^{-1} \ \to \ \Kc M^{-1} L^{-1}|_D \]
where $\hsigma = {^t\sigma} \otimes \Iden_{\Kc M^{-1}}$. The kernel of this is an elementary transformation of $E^*$. As the restriction $\Kc M^{-1} L^{-1} \to \Kc M^{-1} L^{-1}|_D$ is always surjective, the elementary transformation has degree $\deg E^* - e$ if and only if $\sigma$ is a vector bundle injection at each point of $D$. Pulling back to $S^e_M$ as defined in (\ref{SeMdefn}), we obtain a family of elements of $\Quot^e ( E^* )$ parametrised by the open subset
\begin{multline*} \left\{ \left( [ L \xrightarrow{\sigma} \Kc M^{-1} \otimes E ] , D \right) : L \cong \Oc ( -D ) \hbox{ and} \right. \\
 \left. \hbox{$\sigma$ is a vector bundle injection along } D \right\} \end{multline*}
of $S^e_M$. By the universal property of $\Quot^e ( E^* )$, we obtain a map $\geM \colon S^e_M \dashrightarrow \Quot^e ( E^* )$ defined precisely on the above open set, given by sending
\begin{equation} ( \sigma , D ) \ \mapsto \ \Ker \left( E^* \ \xrightarrow{\hsigma} \ \Kc M^{-1} ( D ) \ \to \ \Kc M^{-1} ( D )|_D \right) . \label{gammaDefn} \end{equation}

Now we characterise the locus $\Qeo ( E, M, H^0 ( E^* \otimes M ) ) = \Qeo$ using the spaces $S^e_M$. The description is analogous to that of the osculating spaces to the scroll $\PP E$ given in \cite[Proposition 4.1]{Hit19}.

\begin{proposition} \label{CharSec} Let $E$ and $M$ be as above, and $\geM$ as defined in (\ref{gammaDefn}).
\begin{enumerate}
\item[(a)] We have $\Image ( \geM ) \subseteq \Qeo$.
\item[(b)] Suppose that $[F^* \to E^*]$ belongs to $\Qeo$. Then for some $e' \in \{ 1 ,\ldots , e \}$ there exists $[ G^* \to E^* ]$ belonging to $\Image \left( \gamma^{e'}_M \right) \subseteq \Quot^{e'} ( E^* )$ such that $[ F^* \to G^* ]^* \in \Quot^{e-e_1} ( G^* )$.
\end{enumerate}
\end{proposition}

\begin{proof} (a) Suppose that $[ F^* \to E^* ]$ is of the form $\geM ( \sigma, D)$. 
 By definition of $\geM$, we have a exact diagram
\[ \xymatrix{ F^* \ar[r] \ar[d]^{\hsigma'} & E^* \ar[r] \ar[d]^{\hsigma} & E^* / F^* \ar[d]^\wr \\
 \Kc M^{-1} \ar[r] & \Kc M^{-1} ( D ) \ar[r] & \Kc M^{-1} (D)|_D } \]
Thus $\hsigma'$ defines an element of $H^0 ( \Kc M^{-1} \otimes F )$. Now by hypothesis, $\hsigma$ is a vector bundle surjection along $D$. It follows that $\hsigma'$ does not factorise via $E^*$. Hence $h^0 ( \Kc M^{-1} \otimes F ) \ge h^0 ( \Kc M^{-1} \otimes E ) + 1$. Using (\ref{differentf}), we obtain $h^0 ( F^* \otimes M ) \ge h^0 ( E^* \otimes M ) - e + 1$. Hence $[ F^* \to E^* ]$ belongs to $\Qeo$.

(b) Suppose that $h^0 ( F^* \otimes M ) \ge h^0 ( E^* \otimes M ) - e + 1$. Again using (\ref{differentf}),
 we see that there exists $\sigma \in H^0 ( C, \Kc M^{-1} \otimes F )$ whose image in $H^0 ( C , \Kc M^{-1} \otimes (F/E))$ is nonzero.

Now we view sections of $F$ as rational sections of $E$ with at worst poles limited by $F/E$. (This generalises the definition of $\Oc (D)$ for an effective divisor $D$.) Let $x_1 , \ldots, x_s \in \Supp ( F/E )$ be the points where $( \sigma \mod \Kc M^{-1} \otimes E)$ is supported; as $\sigma \not\in H^0 ( \Kc M^{-1} \otimes E )$, there is at least one such point. For $1 \le j \le s$, we may choose a suitable neighbourhood $U_j$ of $x_j$ and a frame $\nu^{(j)}_1 , \ldots , \nu^{(j)}_r$ for $\Kc M^{-1} \otimes E|_{U_j}$ and a local parameter $z_j$ on $C$ at $x_j$, such that the image of $\sigma|_{U_j}$ is $\cO_{U_j} \cdot \left( z_j^{-e_j} \cdot \nu^{(j)}_1 \right)$. Let $G$ be the elementary transformation of $E$ such that $\Kc M^{-1} \otimes G$ has frame
\begin{equation} z_j^{-e_j} \cdot \nu^{(j)}_1 , \ \nu^{(j)}_2 , \ \ldots , \ \nu^{(j)}_r \label{FrameForG} \end{equation}
over each $U_j$, and is equal to $\Kc M^{-1} \otimes E$ otherwise. Now $z_j^{-e_j} \cdot \nu^{(j)}_1 \subseteq \Image (\sigma)$ belongs to $F$, and the remaining frame elements $\nu^{(j)}_2 , \ldots , \nu^{(j)}_r$ belong to $E \subset F$. Therefore, we have inclusions $E \subset G \subseteq F$. Writing $e' := \sum_{i=1}^s e_j$, we have
\[ \deg E + e' \ = \ \deg G \ \le \ \deg F \ = \ \deg E + e . \]
In particular $e' \le e$, and
\begin{equation} [ F^* \to G^* ] \hbox{ is an element of } \Quot^{e - e'} ( G^* ) . \label{FinFp} \end{equation}

Let $D'$ be the effective divisor $\sum_{i=1}^s e_j x_j \in C_{e'}$. From the construction of $G$ it follows that
\[ \sigma \ \in \ H^0 ( \Kc M^{-1} \otimes G ) \ \subseteq \ H^0 ( \Kc M^{-1} (D') \otimes E ) \]
and moreover, that the image of $\sigma \colon \Oc \to \Kc M^{-1} \otimes G$ is nonzero at each of the $x_j$. Twisting by $\Oc ( -D' )$, we see that $\sigma$ defines a sheaf injection $\sigma' \colon \Oc( -D') \to \Kc M^{-1} \otimes E$ which is a vector bundle injection along $\Supp \, D'$.

Now we claim that
\[ [ G^* \to E^* ] \ = \ \gamma^{e'}_M \left( [ \sigma' \colon \Oc (-D') \to \Kc M^{-1} \otimes E ] , D' \right) . \]
To see this: For each $U_j$, let $\phi_1^{(j)} , \ldots , \phi_r^{(j)}$ be the frame for $E^*|_{U_j}$ dual to $\nu_1^{(j)} , \ldots , \nu_1^{(j)}$. Then by construction of $G$, we see that $G^*|_{U_j} \subset E^*|_{U_j}$ has the frame
\begin{equation} z_j^{e_j} \cdot \phi_1^{(j)} , \ \phi_2^{(j)} , \ \ldots , \ \phi_r^{(j)} \label{frameFp} \end{equation}
dual to (\ref{FrameForG}); and $G^*$ coincides with $E^*$ on $C \setminus D'$. On the other hand, the map
\[ E^* \ \xrightarrow{{^t\sigma'} \otimes \Iden_{\Kc M^{-1}}} \ \Kc M^{-1} ( D' ) \ \to \ \Kc M^{-1} ( D' )|_{D'} \]
is given in our chosen coordinates over $U_j$ by
\[ \sum_{i=1}^r f_i \phi_i^{(j)} \ \mapsto \ z_j^{-e_j} \cdot f_1 \ \mapsto \ z_j^{-e_j} \cdot f_1 \mod \cO_{U_j} . \]
Comparing with (\ref{frameFp}), we see that this map has kernel exactly $G^*$, proving the claim. Statement (b) now follows from the claim and (\ref{FinFp}). \end{proof}

\begin{remark} Suppose $[ L \xrightarrow{\sigma} \Kc \otimes E ]$ is an element of $\qoe$ where $|L^{-1}|$ has dimension $k \ge 1$. For simplicity, suppose $|L^{-1}|$ is base point free. Then $S^e_\Oc$ contains the locus
\[ \{ ( \sigma , D ) : D \in |L^{-1}| \} \ \cong \ \PP^k . \]
By Proposition \ref{CharSec} (a), we obtain an element of $\Qeo$ for each $D \in |L^{-1}|^*$ along which $\sigma$ is a bundle injection. By Lemma \ref{alphaProperties} and Proposition \ref{HQequiv} the projective model $\PP E \dashrightarrow | \cL_M |^*$ has a family of $1$-defective $e$-secants parametrised by an open subset of $|L^{-1}|$. 
\end{remark}

\begin{remark} One may ask whether there is a relation between $S^e_M$ and $\Heo$ akin to that described in Proposition \ref{CharSec} for $\Qeo$. There does exist a map $\beta \colon S^e_M \dashrightarrow \HerS_\nd$ by $( \sigma, D ) \mapsto \PP \sigma (D)$, globalising the construction of Lemma \ref{IdealSection}. This is again defined exactly when $\sigma$ is a vector bundle injection at $D$, and in fact one can show that $\alpha \circ \beta = \gamma$, where $\alpha \colon \Hilb^e ( S ) \dashrightarrow \Quot^e (E^*)$ is the map studied in {\S} \ref{linkHQ}. An analogue of Proposition \ref{CharSec} (a) then follows from Lemma \ref{HQequiv}. However, Proposition \ref{CharSec} (b) is harder to generalise, as one cannot always choose a subscheme playing the role of the sheaf $G$ above; essentially due to issues of identifiability (defined for example in \cite[Definition 6.1]{BC}) arising from the existence of linear spaces contained in $\PP E$. Thus we content ourselves for now with the study of $\Qeo$. \end{remark}

In the sections that follow, we use Proposition \ref{CharSec} to characterise semistability, to study questions of very ampleness of $\opeo$, and to show that the secant loci are well behaved when certain parameters are chosen generally.

\section{The Segre invariant \texorpdfstring{$s_1$}{s1} and secant loci}

In this section, we use the parameter spaces $S^e_M$ to characterise nonempty secant loci for $f = 1$ in terms of the Segre invariant $s_1$. Firstly, we review the notion of Segre invariants.

\subsection{Segre invariants} \label{SegreInvBkgrnd}

\begin{definition} \label{DefnSegreInv} Let $E \to C$ be a vector bundle of rank $r$ and degree $d$. For $1 \le t \le r - 1$, the Segre invariant $s_t (E)$ is defined by
\[ s_t ( E ) \ := \ \min \{ td - r \cdot \deg ( F ) : F \subset E \hbox{ a subbundle of rank } t \} . \]
\end{definition}

Segre invariants are ``degrees of stability'': $E$ is stable if and only if $s_t ( E ) > 0$ for $1 \le t \le r - 1$. It is clear that $s_1 ( E ) = s_1 ( E \otimes L )$ for any line bundle $L$. Moreover, $s_t (E)$ is lower semicontinuous as $E$ varies in families. Taking direct sums of line bundles, one can produce bundles $E$ with arbitrarily low $s_t ( E )$. However, Hirschowitz \cite[Th\'eor\`eme 4.4]{Hir} (see also \cite{CH}) proved the following sharp upper bound on $s_t (E)$.

\begin{theorem}[Hirschowitz's bound] \label{Hirschowitz} Let $E \to C$ be a vector bundle of rank $r$ and degree $d$. Then for $1 \le t \le r - 1$, we have $s_t ( E ) \le t ( n - t )(g-1) + \delta$, where $\delta$ is the unique integer such that $0 \le \delta \leq r-1$ and $t ( r - t ) (g - 1) + \delta \equiv td \mod r$. Moreover, a general stable $E$ attains this upper bound. \end{theorem}

In particular, this implies that Quot schemes of subsheaves of $E$ are always nonempty for low enough degree. Let us make this precise for line subbundles. The following is well known, but we include a proof for convenience.

\begin{lemma} \label{dimQuotOne} Let $C$ be any curve. Let $E \to C$ be a vector bundle which is general in moduli. Then $Q_{1, -e} ( \Kc M^{-1} \otimes E )$ is nonempty and smooth of dimension $re + d + (r+1)(g-1)$ whenever this is nonnegative, and empty otherwise. Moreover, when $\qoe$ is nonempty, a general element is a vector bundle injection. \end{lemma}

\begin{proof} One computes easily that $re + d + (r+1)(g-1) \ge 0$ if and only if
\[ \deg ( \Kc \otimes E ) - r \cdot (-e) \ \ge \ (r-1)(g-1) + \delta , \]
with $\delta$ as in Theorem \ref{Hirschowitz}, using the fact that left side is congruent to $d$ modulo $r$. In this case, by the Hirschowitz bound, $\Kc \otimes E$ has a line subbundle of degree at least $-e$, so $\qoe$ is nonempty. Smoothness and dimension of $\qoe$ for general $E$ then follow from \cite[Lemma 3.3]{LN} and the discussion before it. To see that a general element is saturated, we note that the locus of nonsaturated subsheaves in $\qoe$ has dimension at most
\[ \max\{ \dim Q_{1, -e_1} ( \Kc \otimes E ) + \dim C_{e - e_1} : e_1 < e \} , \]
 which one checks is strictly less than $re + d + (r+1)(g-1)$.

On the other hand, if if $re + d + (r+1)(g-1) < 0$ then $\deg ( \Kc \otimes E ) + re < (r-1)(g-1)$. For such $e$, there is an invertible subsheaf of degree $-e$ in $\Kc \otimes E$ if and only if $s_1 (E)$ is smaller than the generic value, so $E$ is not general. \end{proof}

\subsection{Nonemptiness criterion for \texorpdfstring{$f = 1$}{f = 1}}

The following is a generalisation of \cite[Theorem 5.2]{Hit19}, with a very similar proof.

\begin{theorem} \label{SegreInvChar} Let $E \to C$ be a vector bundle of rank $r$ and degree $d$. For $e \ge 1$, the following are equivalent.
\begin{enumerate}
\renewcommand{\labelenumi}{(\arabic{enumi})}
\item $s_1 (E) > d + r ( 2g - 2 + e )$.
\item For all $M \in \Pic^0 ( C )$ and all $Z \in \HerS_\nd$, the span of $\phi_{\cL_M} ( Z )$ is of dimension $e - 1$ in $| \opeo \otimes \pi^* M |^*$.
\item For all $M \in \Pic^0 ( C )$ and all $[F^* \subset E^*] \in \Quot^e ( E^* )$, we have
\[ h^0 ( F^* \otimes M ) \ = \ h^0 ( E^* \otimes M ) - e . \]
\end{enumerate}
\end{theorem}

\begin{proof} The equivalence of (2) and (3) follows from Proposition \ref{HQequiv} and surjectivity of $\alpha \colon \HerS \to \Quot^e ( E^* )$. 
 Let us show the equivalence of (1) and (3).

Assume (1). As $s_1 ( \Kc M^{-1} \otimes E ) = s_1 ( E ) \ > \ d + r ( 2g - 2 + e )$, if $L \subset \Kc M^{-1} \otimes E$ is an invertible subsheaf, then $\deg ( L ) < -e$. 
 In particular, $Q_{1, -e'} ( \Kc M^{-1} \otimes E )$ and hence $S^{e'}_M$ are empty for all $M \in \Pic^0 (C)$ and $e' \le e$. By Proposition \ref{CharSec} (b), also $\Qeo$ is empty.

Conversely, suppose that $s_1 (E) = s_1 ( \Kc \otimes E ) \le d + r ( 2g - 2 + e )$. Then there exists $L \in \Pic^{-e} (C)$ and a sheaf injection $\sigma' \colon L \to \Kc \otimes E$ of degree $-e$. Let $D \in \Ce$ be any divisor along which $L \to \Kc \otimes E$ is a vector bundle injection. Set $M := L (D) \in \Pic^0 ( C )$ and
\[ \sigma \ := \ \sigma' \otimes \Iden_{L^{-1} ( -D)} \colon \Oc(-D) \ \to \ \Kc M^{-1} \otimes E . \]
Then $( \sigma, D )$ defines an element of $S^e_M$ at which $\gamma^e_M$ is defined. Hence $\Qeo$ is nonempty by Proposition \ref{CharSec} (a).
\end{proof}

\begin{remark} As in \cite[Remark 5.6]{Hit19}, we note that Theorem \ref{SegreInvChar} does not hold for incomplete linear systems. For any $Z \in \HerS$, if we project from a point of $\Span \, \psi (Z)$ then, whatever the value of $s_1 (E)$, the image of $\PP E$ has a defective $e$-secant space. (However, as noted in {\S} \ref{GenProj}, under certain conditions the $\Heo$ and $\Qeo$ do behave as expected under projection from a general centre.) \end{remark}

\noindent We mention now some special cases, generalising those discussed in \cite[Corollary 5.4]{Hit19}. We omit the proofs, as they are practically identical to those in \cite{Hit19}.
Recall that the \textsl{slope} of a bundle $W \to C$ is the ratio $\mu(W) := \deg(W)/\rank(W)$.

\begin{corollary} Let $E$ be a bundle of rank $r$ and degree $d$, where $d \le r ( 1 - 2g )$.
\begin{enumerate}
\item[(a)] Suppose that $s_1 (E) > 0$ (this is the case for example if $E$ is stable). Then for $e \le \mu ( E^* ) - ( 2g - 2 )$, for all $M \in \Pic^0 (C)$ the loci $\left( \Heo \right)_\nd$ and $\Qeo$ are empty.
\item[(b)] If $e \ge \mu(E^*) - \frac{1}{r}(r+1)(g-1)$, then for some $M \in \Pic^0 (C)$, the loci $\left( \Heo \right)_\nd$ and $\Qeo$ are nonempty.
\item[(c)] If $s_1 (E)$ is the generic value $(r-1)(g-1) + \delta$, then the converse to (b) also holds.
\end{enumerate}
\end{corollary}

\subsection{A question of Lange} \label{Lange} If $r = 2$, then $\PP E$ is a ruled surface over $C$. For rational curves and elliptic curves, Lange \cite[Lecture 2 (a)]{Lan} (see also \cite[Exercise V.2.12]{Har}) gives criteria for very ampleness of $\opeo \otimes \pi^* M$ in terms of $\deg (M)$ and $s_1 (E)$. He then poses the problem of finding similar criteria for ruled surfaces when $g \ge 2$. We now offer a result in this direction, which in fact holds for projective bundles of any dimension $r \ge 2$ over $C$.

\begin{theorem} \label{LangeQ} Let $E \to C$ be a bundle of rank $r$ and degree $d$. Then the line bundle $\opeo \otimes \pi^* L$ is very ample for all $L \in \Pic^\ell ( C )$ if and only if
\begin{equation} \ell \ > \ \mu( E ) - \frac{s_1 ( E )}{r} + 2g . \label{soEmb} \end{equation}
\end{theorem}

\begin{proof} A map $S \to \PP^n$ is an embedding if and only if all $2$-secants are nondefective. Furthermore, when $S$ is a projective bundle $\PP E$ over $C$, it follows from the argument of \cite[Theorem 3 (i)]{StaQS} that all $2$-secants are $\pi$-nondefective (cf.\ Definition \ref{DefnPiDef}).
 Thus it suffices to show that (\ref{soEmb}) is equivalent to emptiness of the secant locus $H^1_2 ( \opeo \otimes \pi^* L )_\nd$ for all $L \in \Pic^\ell ( C )$.

Fix now a line bundle $L_0$ of degree $\ell$. For any $M \in \Pic^ (C)$, we have
\[ H^1_2 ( H^0 ( \PP E , \opeo \otimes \pi^* L_0 M ) ) \ \cong \ H^1_2 \left( H^0 ( \PP ( E \otimes L_0^{-1} ) , \cO_{\PP ( E \otimes L_0^{-1} )} ( 1 ) \otimes \pi^* M )) \right) . \]
By Theorem \ref{SegreInvChar}, the right-hand side is empty for all $M \in \Pic^0 (C)$ if and only if
\[ s_1 ( E \otimes L_0^{-1} ) \ > \ d - r\ell + 2rg . \]
As $s_1 ( E \otimes L_0^{-1} ) = s_1 ( E )$, this is equivalent to 
 $\ell > \mu ( E ) - \frac{s_1 (E)}{r} + 2g$, as desired. \end{proof}

Theorem \ref{LangeQ} is a statement concerning $E \otimes L$ for all line bundles $L$ of a given degree. We make one more observation showing that very ampleness may also depend on $L \in \Pic^\ell (C)$. Fix such an $L$, and write $h := h^0 ( \Kc L^{-1} \otimes E )$. We consider the twisted Brill--Noether locus
\[ B^{h+1}_{1, 2} ( \Kc L^{-1} \otimes E ) \ = \ \{ N \in \Pic^2 (C) : h^0 ( \Kc L^{-1} N \otimes E ) \ge h + 1 \} . \]
Also on $\Pic^2 (C)$ we have the standard Brill--Noether locus
\[ B^1_{1, 2} \ = \ \{ \Oc ( x + y ) : x , y \in C \} \ \subseteq \ \Pic^2 ( C ) \]
parametrising effective line bundles of degree two over $C$.

\begin{proposition} Let $E$ be a bundle of rank $r$ and degree $d$. For $L \in \Pic^\ell ( C )$, the line bundle $\opeo \otimes \pi^* L$ is very ample if and only if
\[ B^{h+1}_{1, 2} ( \Kc M^{-1} \otimes E ) \cap B^1_{1, 2} \ = \ \emptyset . \]
\end{proposition}

\begin{proof} By \cite[Proposition 4.1]{Hit20} (substituting $\Kc M^{-1} \otimes E$ for ``$V$''), the line bundle $\opeo \otimes \pi^* L$ is very ample on $\PP E$ if and only if
\begin{equation} h^0 ( C, \Kc L^{-1} \otimes E (x+y) ) \ = \ h^0 ( C, \Kc L^{-1} \otimes E ) \hbox{ for all } x + y \in C_2 . \label{CondOne} \end{equation}
This is equivalent to saying that $h^0 ( C, \Kc L^{-1} \otimes E (x+y) ) = h$ for all $x + y \in C_2$; in other words, that no line bundle of the form $\Oc ( x + y )$ belongs to the locus $B^{h + 1} ( \Kc L^{-1} \otimes E )$. \end{proof}

\subsection{A criterion for semistability} \label{SemistCohomSt}

Here we use Theorem \ref{SegreInvChar} to give a criterion for semistability. This is entirely analogous to \cite[Theorem 5.7]{Hit19}; as the proof is practically identical, we only sketch it. For a given $E \to C$, for $t \ge 1$ we write $\pi_t \colon \PP ( \wedge^t E ) \to C$ for the projection.

\begin{theorem} Let $E$ be a vector bundle of rank $r$ with $\mu (E) < 1 - 2g$. Then the following are equivalent.
\begin{enumerate}
\item[(1)] The bundle $E$ is semistable.
\item[(2)] For $1 \le t \le r-1$, for all $M \in \Pic^0 (C)$, and for $1 \le e < t \cdot \mu( E^* ) - (2g-2)$, the secant locus
\[ H^{e-1}_e \left( H^0 ( \PP ( \wedge^t E ) , \cO_{\PP ( \wedge^t E )} (1) \otimes \pi_t^* M \right )_\nd \]
is empty.
\item[(3)] For $1 \le t \le r-1$, and for $M$ and $e$ as in (2), 
 the secant locus
\[ Q^{e-1}_e \left( \wedge^t E , M , H^0 ( C, \wedge^t E \otimes \pi_t^* M ) \right) \]
is empty; that is, $h^0 ( C, \cF_n ) = h^0 ( C, M \otimes \wedge^n E^* ) - e$ for all elementary transformations $0 \to \cF_n \to \wedge^n E^* \to T \to 0$ where $T$ is torsion of length $e$.
\end{enumerate}
\end{theorem}

\begin{proof} The equivalence of (2) and (3) is proven as in Theorem \ref{SegreInvChar}. If $E$ is semistable then, since $\C$ has characteristic zero, also $\wedge^t E$ is semistable for $1 \le t \le r-1$. In particular, $s_1 ( \wedge^t E ) \ge 0$. Theorem \ref{SegreInvChar} then implies (3).

Conversely, assume (3). Then for each $t \in \{ 1, \ldots , r-1 \}$, by Theorem \ref{SegreInvChar} we have
\[ s_1 ( \wedge^t E ) \ > \ \deg ( \wedge^t E ) + \rank ( \wedge^t E ) \cdot (2g - 2 - e_t) \]
where $e_t \ = \ \max\{ e : 1 \le e < t \cdot \mu ( E^* ) - (2g - 1) \}$. By definition of $s_1$, the left and right sides of the last inequality are congruent modulo $\rank ( \wedge^t E )$. Therefore
\[ s_1 ( \wedge^t E ) \ \ge \ \deg ( \wedge^t E ) + \rank ( \wedge^t E ) \cdot ( 2g - 2 + e_t + 1 ) . \]
By definition of $e_t$, this becomes $s_1 ( \wedge^t E ) \ge 0$. 
 As in the proof of \cite[Theorem 5.7]{Hit19}, we conclude that $\mu ( F ) \le \mu ( E )$ for all rank $t$ subbundles $F \subset E$. \end{proof}

\section{Secant loci of general scrolls for \texorpdfstring{$f = 1$}{f = 1}} \label{General}

In this final section we study questions of expected dimension, nonemptiness and enumeration for the secant loci $\Qeo ( E, M, V ) =: \QeoV$ when the parameters $E$, $M$ and $V$ are chosen generally. The statements are valid for any curve $C$. We obtain also the emptiness of certain $\Heo ( \cL_M , V )_\nd$; however, in the lack of a complete description of the fibres of $\alpha \colon \Hilb^e ( \PP E ) \dashrightarrow \Quot^e ( E^* )$, we do not prove further dimension bounds on $\Heo ( \cL_M, V )$ at present.

\subsection{Expected dimension of \texorpdfstring{$\QeoV$}{Q e-1 e (V)}}

Here we show that if the parameters are chosen generally, $\left( \HeoV \right)_\nd$ and $\QeoV$ are empty if the expected dimension is zero, and $\QeoV$, if nonempty, is of the expected dimension otherwise. This generalises \cite[Theorem 6.4]{Hit19}, and uses a similar method. We begin by treating the complete case. Fix $r \ge 2$ and $d < r( 1 - g ) - 1$, so that $n := -d - r(g-1) - 1 \ge 0$. Recall that for $f = 1$, the expected dimension of $\Heo$ and $\Qeo$ is $re - ( n + 1 - e + 1) = (r+1)e - n - 2$.

\begin{theorem} \label{GeneralGood} Let $E \to C$ be a general vector bundle of rank $r$ and degree $d$. Then there is a nonempty open subset $U \in \Pic^0 (C)$ such that for $M \in U$, the following hold.
\begin{enumerate}
\item[(a)] $h^0 ( C, E^* \otimes M ) = n + 1$, so $| \cL_M |$ has dimension $n$.
\item[(b)] If $(r+1)e - n - 2 < 0$, then $\Qeo$ and $\left( \Heo \right)_\nd$ are empty.
\item[(c)] If $(r+1)e - n - 2 \ge 0$, then $\Qeo$ is empty or of the expected dimension.
\end{enumerate} \end{theorem}

\begin{proof} Part (a) is exactly \cite[Theorem 6.4 (a)]{Hit19}. For (b): We recall the construction of $S^e_M$ as a fibre product in (\ref{SeMdefn}), and the forgetful map
\[ b \colon Q_{1, -e} ( \Kc M^{-1} \otimes E ) \ \to \ \Pic^{-e} (C) . \]
Unwinding definitions, we see for any $M \in \Pic^0 (C)$ that
\[ M^{-1} \cdot b \left( Q_{1, -e} ( \Kc \otimes E ) \right) \ = \ b \left( Q_{1, -e} ( \Kc M^{-1} \otimes E ) \right) , \]
where the action on the left hand side is that of $\Pic^0 ( C )$ on itself by translation. Thus, by Theorem \ref{KleimanTheorem}, for general $M \in \Pic^0 ( C )$ the space $S^e_M$ is empty or of dimension
\[ \dim Q_{1, -e} ( \Kc M^{-1} \otimes E ) + \dim \Ce - \dim \Pic^{-e} (C) \]
when this is nonnegative. By Lemma \ref{dimQuotOne} and the assumption of generality of $E$, this number is
\[ re + d + (r+1)(g-1) + e - g \ = \ 
 (r+1)e - n - 2 , \]
the last equality by the definition of $n$ above. This is exactly the expected dimension of $\Qeo$ and $\Heo$. 

Now by Proposition \ref{CharSec} (b), the locus $\Qeo$ is nonempty only if $S^{e'}_M$ is nonempty for some $e' \le e$. The above computation shows that $\dim S^{e'}_M < \dim S^e_M$ for $e' < e$. Therefore, if $(r+1)e - n - 2 < 0$ then $\Qeo$ is empty. As $\alpha \left( \Heo \right)_\nd = \Qeo$ by Proposition \ref{HQequiv}, also $\left( \Heo \right)_\nd$ is empty. This proves (b).

For the rest: If $[ F^* \to E^* ]$ is a point of $\Qeo$ then, by Proposition \ref{CharSec} (b), we have $F^* \subseteq G^*$, where $[ G^* \to E^* ]$ belongs to the image of
\[ \gamma^{e'}_M \colon S^{e'}_M \ \dashrightarrow \ \Quot^{e'} (E^*) , \]
for some $e' \le e$. Thus $\dim \Qeo$ is at most
\[ \dim S^{e'}_M + \dim \Quot^{e-e'} ( G^* ) \ = \ (r+1)e' -n-2 + r (e-e') \ = \ (r+1)e - n - 2 - (e - e') . \]
As $e \ge e'$, part (c) follows. (Note that we have equality for $e' = e$; compare with Proposition \ref{firstprops} (b).) \end{proof}

Using Corollary \ref{QefGenProj}, we can extend Theorem \ref{GeneralGood} to the case of incomplete linear systems.

\begin{theorem} \label{GeneralGoodIncomplete} Let $E$ and $M$ be general in the sense of Theorem \ref{GeneralGood}. Let $W \subseteq H^0 ( C, E^* \otimes M )$ be a general subspace of dimension $m + 1$.
\begin{enumerate}
\item[(a)] If $(r+1)e - m - 2 < 0$, then $\Qeo (W)$ is empty.
\item[(b)] If $(r+1)e - m - 2 < 0$, then $\left( \Heo (W) \right)_\nd$ is empty.
\item[(c)] If $(r+1)e - m - 2 \ge 0$, then $\Qeo (W)$ is empty or of the expected dimension.
\end{enumerate} \end{theorem}

\begin{proof} For (a) and (c), we will apply Corollary \ref{QefGenProj}. As $f = 1$, to satisfy the hypothesis it will suffice to show that $\Qeo ( H^0 ( E^* \otimes M ) )$ is empty or has the expected dimension (trivially, $Q^e_e$ is always of the expected dimension). This follows from Theorem \ref{GeneralGood}. Then (b) follows as before from Lemma \ref{alphaProperties} (a) and Proposition \ref{HQequiv}. \end{proof}

\subsection{Nonemptiness of \texorpdfstring{$\QeoV$}{Q e-1 e (V)}} \label{nonemptiness}

Theorems \ref{GeneralGood} and \ref{GeneralGoodIncomplete} show that $\QeoV$, if nonempty, is of the expected dimension for a general choice of the parameters. We now describe two situations in which $\QeoV$ is in fact nonempty.

\begin{proposition} Let $E \to C$ be a bundle which is general in moduli. Suppose $e \ge \max\{ g , -d/r - (g-1) (r+1)/r \}$. Then $\QeoV$ is nonempty and of the expected dimension for generic $M \in \Pic^0 ( C )$. \end{proposition}

\begin{proof} As $e \ge -d/r - (g-1) (r+1)/r$, we have $re + d + (r+1)(g-1) \ge 0$. Since $E$ is general, by Lemma \ref{dimQuotOne} the scheme $\qoe$ is nonempty and a general element is of the form $[ L \xrightarrow{\sigma} \Kc M^{-1} \otimes E ]$ where $\sigma$ is a vector bundle injection. But since also $e \ge g$, for any $L$ so occurring we have $L = \Oc ( -D )$ for some $D \in \Ce$. Therefore, $( \sigma , D )$ is a point of $S^e_M$ where $\geM$ is defined. Hence $\Qeo$ is nonempty by Proposition \ref{CharSec} (a). \end{proof}

\begin{proposition} \label{NonemptinessTwo} Let $E$ be any bundle, and suppose that $\qoe$ is nonempty for some $e \ge 1$. Then $\Qeo ( E, M, H^0 ( E^* \otimes M ) )$ is nonempty for some $M \in \Pic^0 (C)$. \end{proposition}

\begin{proof} Suppose $[ L \xrightarrow{\sigma} \Kc \otimes E ]$ is a point of $\qoe$. Now $\sigma$ may not be a bundle injection at all points, but we can choose a sufficiently general $D \in \Ce$ along which $\sigma$ is a bundle injection. Then $L^{-1} ( -D )$ has degree zero, and
\[ \left( \left[ \Oc ( -D ) \xrightarrow{\sigma \otimes L^{-1} ( -D )} \Kc L^{-1} (-D) ) \otimes E \right] , D \right) \]
is a point of $S^e_{L^{-1} ( -D )}$ at which $\gamma^e_{L^{-1} ( -D)}$ is defined. Then
\[ \Qeo ( E , L^{-1} (-D) , H^0 ( C, E^* \otimes L^{-1}(-D) ) ) \]
is nonempty by Proposition \ref{CharSec} (a). \end{proof}

\begin{remark} Proposition \ref{NonemptinessTwo} applies whenever $\qoe$ has nonnegative expected dimension. If $E$ is generic, then this follows from Lemma \ref{dimQuotOne}. For arbitrary $E$, one can find possibly nonsaturated elements of $\qoe$ by taking subsheaves of $Q_{1, h} ( \Kc \otimes E )$ for $h > -e$. \end{remark}

\subsection{Enumeration of \texorpdfstring{$\QeoV$}{Q e-1 e}}

By the previous two subsections, one knows that $\Qeo ( E, M, V ) = \QeoV$, if nonempty, is of the expected dimension for any $C$ and for a general choice of $(E, M, V)$; and furthermore that $\QeoV$ is indeed nonempty under certain conditions. We now use results of \cite{EGL}, \cite{OP} and \cite{StaCL} to enumerate $\QeoV$ when it has and attains expected dimension zero.

Following the notation of \cite{OP} and \cite{StaCL}, we denote by $\Me$ the vector bundle over $\Quot^e ( E^* )$ with fibre $H^0 ( C, (E^*/ F^*) \otimes M )$ at the point $[F^* \to E^*]$. For $V \subseteq H^0 ( E^* \otimes M )$, let $\varepsilon \colon \cO_{\Quot^e (E^*)} \otimes V \to \Me$ be the evaluation map. Then $\QeoV$ is the determinantal locus $\{ F^* \in \Quot^e ( E^* ) : \rank \, \varepsilon|_{F^*} \le e - 1 \}$. When $re = n + 2 - e$, so that $\QeoV$ has expected dimension zero, using the Porteous formula we obtain
\begin{equation} [ \QeoV ] \ = \ s_{n+2-e} \left( \Me \right) \ = \ \int_{[ \Quot^e ( E^* ) ]} s ( \Me ) . \label{integral} \end{equation}
The main goal of this subsection will be to sketch the proof of the following.

\begin{theorem} \label{counting} Let $C$ be a smooth curve, $E \to C$ a vector bundle of rank $r$ and $M \to C$ a line bundle. Then there is an equality of formal power series
\begin{equation} \sum_{e \ge 0} \int_{[ \Quot^e_C ( E^* ) ]} s ( \Me ) \, q^e \ \ = \ A_1 ( q )^{\deg (E^* \otimes M)} \cdot B ( q )^{1-g} \label{GenSeries} \end{equation}
where
\[ q = (-1)^r t ( 1 + t )^r , \quad A_1 ( q ) = 1 + t \quad \hbox{and} \quad B ( q ) = \frac{(1 + t)^{r+1}}{1 + t ( r+1 )} . \]
\end{theorem}

\begin{proof}[Sketch of proof] All the ingredients for proving this theorem are present in \cite{EGL}, \cite{OP} and \cite{StaCL}, so we give only an outline. The formula (\ref{GenSeries}) is proven in \cite[Theorem 6 and Theorem 8]{OP} when $E$ is the trivial bundle of rank $r \ge 1$. The generalisation of the corresponding formula \cite[Corollary 16]{OP} for surfaces to arbitrary $E$ is \cite[Proposition 4.3]{StaCL}, and we use the same strategy.

To show the existence of the power series $A_1$ and $B$, we follow essentially word for word the proof of \cite[Theorem 4.2]{EGL}, replacing $\cK_r$ with the set
\[ \cT_r \ := \ \{ ( C, E, M ) : \hbox{$C$ a smooth curve, $E \to C$ a bundle of rank $r$ and $M \in \Pic ( C )$} \} , \]
and $H_{\Psi , \Phi} ( S, x )$ with the map $\cT_r \to \bQ [[q]]$ given by
\[ ( C , E , M ) \ \mapsto \ \sum_{e \ge 0} q^e \int_{[ \Quot^e_C ( E^* ) ]} s ( \Me ) ; \]
and $\gamma$ with the map $(C, E, M) \mapsto ( \deg (E^* \otimes M) , \chi ( \Oc ) )$. The role of \cite[Theorem 4.1]{EGL} is played in our situation by the following statement which emerges in the proof of \cite[Theorem 3.3]{StaCL}:

\begin{proposition} \label{UnivPoly} For each $e \ge 0$, there exists a polynomial $u_e ( x, y )$ which is universal in the sense that for any $(C, E, M) \in \cT_r$ we have
\[ \int_{[\Quot^e (E^*)]} s ( \Me ) \ = \ u_e \left( \deg ( E^* \otimes M ) ,  \chi ( \Oc ) \right) . \]
\end{proposition}

Once the existence of $A_1$ and $B$ is established, they can be computed by comparing with \cite[Theorem 6 and Theorem 8]{OP} for suitable triples $(C, E, M)$ with $E$ trivial and $C$ rational or elliptic. \end{proof}

The following enumeration result is immediate from Theorem \ref{counting} and (\ref{integral}).

\begin{corollary} Suppose $E$, $M$ and $V \subseteq H^0 ( E^* \otimes M )$ are such that $\QeoV$ has and attains expected dimension zero. Then the number of points of $\QeoV$, counted with multiplicity, is the coefficient of $q^e$ in the expression (\ref{GenSeries}). \end{corollary}

\end{document}